\documentclass[aop,MSNbibl,seceqn,citesort,dvips]{arximspdf}

\doi{10.1214/11-AOP732} 
\volume{41}
\issue{1}
\pubyear{2013}
\firstpage{206}
\lastpage{228}

\makeatletter

\newcommand{\cal}{\mathcal}
\newcommand{\mod}{\operatorname{mod}}

\newtheorem{Thm}{Theorem}[section]
\newtheorem{Lem}[Thm]{Lemma}
\newtheorem{Prop}[Thm]{Proposition}
\newtheorem{Coro}[Thm]{Corollary}

\newproclaim{Rem}[Thm]{Remark}
\newproclaim{Def}[Thm]{Definition}
\newproclaim{Example}[Thm]{Example}

\newcommand{\ind}{{\mathbf1}}
\newcommand{\proba}{\mathbb P}
\newcommand{\esp}{{\mathbb E}}
\newcommand{\supp}{\operatorname{supp}}
\newcommand{\defe}{:=}
\newcommand{\inv}{^{-1}}

\newcommand{\calA}{{\cal A}}
\newcommand{\calB}{{\cal B}}
\newcommand{\filF}{{\cal F}}
\newcommand{\calG}{{\cal G}}
\newcommand{\calI}{{\cal I}}
\newcommand{\calT}{{\cal T}}

\newcommand{\lone}{{L^1}}

\newcommand{\la}{L^\alpha}

\newcommand{\sumn}{\sum_{n=1}^\infty}
\newcommand{\sumjp}{\sum_{j=1}^p}
\newcommand{\bcupn}{\bigcup_{n=1}^\infty}

\newcommand{\eqd}{\stackrel{\mathrm{d}}{=}}
\newcommand{\Eintt}{\int^{\!\!\!\!\!\!\!e}}
\newcommand{\Einttt}{\,\,\int^{\!\!\!\!\!\!\!e}\,\,}

\newcommand{\lap}{L^\alpha_+}

\newcommand{\vv}[1]{{\mathbf #1}}

\newcommand{\mam}{M_{\alpha,\vee}}
\renewcommand{\d}{{\mathrm d}}
\newcommand{\bbN}{{\mathbb N}}

\newcommand{\sbm}{{(S,\calB,\mu)}}
\newcommand{\smu}{{(S,\mu)}}

\newcommand{\limn}{\lim_{n\to\infty}}
\newcommand{\limT}{\lim_{T\to\infty}}

\newcommand{\btt}{B(T)}
\newcommand{\ctt}{C(T)}

\newcommand{\bX}{\mathbf{X}}
\newcommand{\Zd}{\mathbb{Z}^d}
\newcommand{\rpd}{\mathbb{R}_{\scriptscriptstyle{+}}^d}
\newcommand{\zpd}{\mathbb{Z}_{\scriptscriptstyle{+}}^d}

\newcommand{\what}{\widehat}
\newcommand{\wtilde}{\widetilde}

\makeatother

\begin{document}
\begin{frontmatter}

\title{Ergodic properties of sum- and max-stable stationary random
fields via null and positive group actions}
\runtitle{Ergodic properties of sum- and max-stable fields}

\begin{aug}
\author[A]{\fnms{Yizao} \snm{Wang}\thanksref{t1}\ead[label=e1]{yizwang@umich.edu}},
\author[B]{\fnms{Parthanil} \snm{Roy}\corref{}\ead[label=e2]{parthanil@isical.ac.in}}
\and
\author[A]{\fnms{Stilian A.} \snm{Stoev}\thanksref{t1}\ead[label=e3]{sstoev@umich.edu}}
\runauthor{Y. Wang, P. Roy and S. A. Stoev}
\affiliation{University of Michigan, Indian Statistical Institute,
Kolkata,
and~University~of~Michigan}
\address[A]{Y. Wang\\
S. A. Stoev\\
Department of Statistics\\
University of Michigan\\
Ann Arbor, Michigan 48109-1107\\
USA\\
\printead{e1}\\
\phantom{E-mail: }\printead*{e3}}
\address[B]{P. Roy\\
Statistics and Mathematics Unit\\
Indian Statistical Institute\\
203 B. T. Road\\
Kolkata 700108\\
India\\
\printead{e2}} 
\end{aug}

\thankstext{t1}{Supported in part by NSF Grants DMS-08-06094 and
DMS-11-06695 at the University of Michigan.}

\received{\smonth{10} \syear{2009}}
\revised{\smonth{12} \syear{2010}}

%
\begin{abstract}
We establish characterization results for the ergodicity of stationary
symmetric $\alpha$-stable (S$\alpha$S) and $\alpha$-Fr\'echet random
fields. We show that the result of Samorodnitsky [\textit{Ann. Probab.}
\textbf{33} (2005) 1782--1803] remains valid in the multiparameter
setting, that is, a stationary S$\alpha$S ($0<\alpha<2$) random field
is ergodic (or, equivalently, weakly mixing) if and only if it is
generated by a null group action. Similar results are also established
for max-stable random fields. The key ingredient is the adaption of a
characterization of positive/null recurrence of group actions by
Takahashi [\textit{K\=odai Math. Sem. Rep.} \textbf{23} (1971)
131--143], which is dimension-free and different from the one used by
Samorodnitsky.
\end{abstract}

%

\end{frontmatter}

\section{Introduction}

A process is called sum-stable (max-stable, resp.) if so are its
finite-dimensional
distributions and it arises as a limit,
under suitable affine transformations, of sums (maxima, resp.)
of independent processes. Convenient stochastic integral
representations have been developed and actively used to study the
structure and properties of sum-stable processes and random fields
(see, e.g., Samorodnitsky and Taqqu~\cite{samorodnitsky94stable},
Rosi\'nski~\cite{rosinski95structure,rosinski00decomposition}, Rosi\'nski
and Samorodnitsky~\cite{rosinski96classes}, Pipiras and Taqqu \cite
{pipiras04stable}, Samorodnitsky \cite
{samorodnitsky04extreme,samorodnitsky04maxima,samorodnitsky05null}, Roy
and Samorodnitsky~\cite{roy08stationary} and Roy \cite
{roy10ergodic,roy10nonsingular}).
On the other hand, the seminal works of de Haan~\cite{dehaan84spectral}
and de Haan and Pickands~\cite{dehaan86stationary} as well as the recent
developments by Stoev and Taqqu~\cite{stoev06extremal}, Wang and
Stoev~\cite{wang10association,wang10structure} and Kabluchko \cite
{kabluchko09spectral} have developed similar tools to represent
and handle general classes of max-stable processes.

The ergodic properties of stationary stochastic processes and random
fields are of fundamental importance and hence well-studied. See, for
example, Maruyama~\cite{maruyama70infinitely}, Rosi\'nski and \.Zak
\cite{rosinski96simple,rosinski97equivalence} and Roy \cite
{roy07ergodic,roy09poisson} for results on infinite divisible processes
and Cambanis et al.~\cite{cambanis87ergodic}, Podg\'orski
\cite{podgorski92note}, Gross and Robertson~\cite{gross93ergodic} and
Gross~\cite{gross94some} for results on stable processes. These
culminated in the characterization of Samorodnitsky \cite
{samorodnitsky05null}, which shows that the ergodicity of a stationary
symmetric stable process is equivalent to the null-recurrence of the
underlying nonsingular flow. On the other hand, the ergodic properties
of max-stable processes have been recently studied by Stoev
\cite{stoev08ergodicity}, Kabluchko~\cite{kabluchko09spectral} and
Kabluchko and Schlather~\cite{kabluchko10ergodic}. In particular,
Kabluchko~\cite{kabluchko09spectral} has shown that as in the
sum-stable case, one can associate a nonsingular flow to the stationary
max-stable process and that the characterization of Samorodnitsky
\cite{samorodnitsky05null} remains valid. The case of random fields,
however, remained open in both sum- and max-stable settings.

Our goal in this paper is to establish a Samorodnitsky-type
characterization for sum-stable and max-stable \textit{random fields}.
The main obstacle is the unavailability of a higher-dimensional
analogue of the work of Krengel~\cite{krengel67classification}, which
plays a crucial role in Samorodnitsky's approach for processes.
We resolve this problem by providing an alternative dimension-free
characterization of ergodicity for both classes of sum- and max-stable
stationary random fields.
For simplicity of exposition as well as mathematical tractability, we
work with symmetric $\alpha$-stable (S$\alpha$S), ($0<\alpha<2$)
sum-stable random fields and $\alpha$-Fr\'echet max-stable random
fields ($\alpha>0$).

The key ingredient of our results is the adaptation of the work of
Takahashi~\cite{takahashi71invariant}. Thanks to Takahashi's result, we
are able to develop tractable
and dimension-free criteria for verifying whether a given spectral
representation corresponds to an S$\alpha$S random field generated by a
null (or positive) action.
We also extend a well-known result of Gross~\cite{gross94some} and give
necessary and sufficient condition for a stationary S$\alpha$S random
field to be weakly mixing and in the process fill a gap in the proof of
\cite{gross94some} (see Remark~\ref{remarkmistakegross} below).
Similar results for $\alpha$-Fr\'echet random fields are obtained.
Furthermore, these results offer alternative characterizations of
ergodicity in the one-dimensional case.

The paper is organized as follows. In Section~\ref{secprelim} we start
with some auxiliary results from ergodic theory. In Section \ref
{secdecomp} we establish the positive-null decomposition for measurable
stationary S$\alpha$S random fields. Section~\ref{secergo}
characterizes the ergodicity of S$\alpha$S random fields. The
max-stable setting is discussed in Section~\ref{secmax-stable}. We
conclude with a couple of examples in Section~\ref{secexample}. Some
technical proofs and auxiliary results are given in the
\hyperref[app]{Appendix}.


\section{Preliminaries on ergodic theory}\label{secprelim}

Throughout this paper, we let $\sbm$ denote a standard Lebesgue space
(see Appendix A in~\cite{pipiras04stable}).
Let $\phi$ denote a bi-measurable and invertible transformation on $S$.
We say that $\phi$ is
\textit{nonsingular}, if the measure $\mu\circ\phi\inv$ and $\mu$ are
equivalent, written $\mu\circ\phi\inv\sim\mu$. In this case, one can
define the \textit{dual operator}
$\widehat\phi$ as a mapping from $L^1(S,\mu)$ to $L^1(S,\mu)$:
%
\begin{equation}\label{eqdual}
\widehat\phi f(s)\equiv[\widehat\phi f](s) \defe
\biggl(\frac{\d(\mu\circ\phi\inv)}{\d\mu}\biggr)(s)
f\circ\phi\inv(s).
\end{equation}
Note that $\widehat\phi$ is a positive linear isometry (hence a
contraction) on $\lone\smu$. The
characterization results in the next section are in terms of dual operators.

\subsection{Group actions}

Let $G\equiv(G,+)$ be a locally compact, topological Abelian group
with identity element $0$. Equip $G$
with the Borel $\sigma$-algebra~$\calA$.

\begin{Def}\label{defgroupAction} A collection of measurable
transformations $\phi_t\dvtx S\to S$, \mbox{$t\in G$} is called a \textit{group
action} of $G$ on $S$ (or a \textit{$G$-action}), if:
\begin{longlist}[(iii)]
\item[(i)] $\phi_0(s) = s$ for all $s\in S$,
\item[(ii)] $\phi_{v+u}(s) = \phi_u\circ\phi_v(s) \mbox{ for all
} s\in
S, u,v\in G$,
\item[(iii)] $(s,u)\mapsto\phi_u(s)$ is measurable w.r.t. the product
$\sigma$-algebra $\calA\otimes\calB$.
\end{longlist}
\end{Def}

A $G$-action $\calG= \{\phi_t\}_{t\in G}$ on $(S,\mu)$ is
\textit{nonsingular} if $\phi_t$ is nonsingular for all $t\in G$.
In this paper, all the group actions are assumed to be nonsingular.

The existence of
a $\calG$-invariant finite measure $\nu$, $\nu\sim\mu$ (equivalently,
the existence of a fixed point
of the dual operator $\widehat\phi$, see, e.g., Proposition
1.4.1 in~\cite{aaronson97introduction}),
is an important problem in ergodic theory.
The investigation of this problem was initiated by Neveu
\cite{neveu67existence} and further explored by Krengel
\cite{krengel67classification} and Takahashi~\cite{takahashi71invariant},
among others. In the rest of this section we present results
due essentially to Takahashi~\cite{takahashi71invariant}. We will see
that the invariant finite measures induce a modulo $\mu$ unique
decomposition of $S$. This decomposition will play an important role in
the characterization of
ergodicity for sum- and max-stable random fields. The proofs of the
results mentioned in this section are given in the \hyperref[app]{Appendix}.

Consider the class of finite (positive) $\calG$-invariant measures on
$S$ absolutely continuous with
respect to $\mu$:
\[
\Lambda(\calG)\defe\{\nu\ll\mu\dvtx \nu\mbox{ finite measure on }
S, \nu
\circ\phi\inv= \nu\mbox{ for all }\phi\in\calG\}.
\]
For all $\nu\in\Lambda(\calG)$, let $S_\nu\equiv\supp(\nu):=\{
\d\nu
/\d\mu>0\}$ denote the support of $\nu$
($\mathrm{mod}\ \mu$) and set $I(\calG)\defe\{S_\nu\dvtx\nu\in\Lambda
(\calG)\}$.
%
\begin{Lem}\label{lemmaximal}
There exists a modulo $\mu$ unique maximal element $P_\calG\in
I(\calG
)$, that is:
\begin{longlist}[(ii)]
\item[(i)] For all $S_\nu\in I(\calG)$, $S_\nu\subset P_\calG$, that
is, $\mu(S_\nu\setminus P_\calG) = 0$.
\item[(ii)] If \textup{(i)} holds for $Q_\calG\in I(\calG)$, then $P_\calG=
Q_\calG$ modulo $\mu$.\vadjust{\goodbreak}
\end{longlist}
\end{Lem}

This result suggests the decomposition
%
\begin{equation}\label{ePos-Null}
S = P_{\calG}\cup N_{\calG},
\end{equation}
where $N_{\calG}:= S\setminus P_{\calG}$. The set $P_\calG\equiv
S_{\nu
_0}, \nu_0\in\Lambda(\calG)$ is the largest
(mod $\mu$) set where one can have a finite $\calG$-invariant measure
$\nu_0$, equivalent to $\mu\vert_{P_\calG}$.
Consequently, there are no finite measures supported on $N_{\calG}$,
invariant w.r.t. $\calG$ and absolutely continuous w.r.t. $\mu$.
The next theorem provides a convenient characterization of the
decomposition (\ref{ePos-Null}).
%
\begin{Thm}\label{thmdecompPN}
Consider any $f\in\lone\smu, f>0$. Let $P_\calG$ denote the unique
maximal element of $I(\calG)$ and set $N_\calG\defe S\setminus
P_\calG
$. We have the following:
\begin{longlist}[(iii)]
\item[(i)] The sets $P_\calG$ and $N_\calG$ are invariant w.r.t.
$\calG$, that is, for all $\phi\in\calG$,
we have
\[
\mu\bigl(\phi\inv(P_\calG)\triangle P_\calG\bigr) = 0
\quad\mbox{and}\quad\mu\bigl(\phi\inv (N_\calG )\triangle N_\calG\bigr) = 0.
\]
\item[(ii)] Restricted to $P_\calG$,
%
\begin{equation}\label{eqPG}
\sumn\widehat\phi_{u_n}f(s) = \infty,\qquad\mu\mbox{-a.e.}\mbox{ for all } \{
\phi
_{u_n}\}_{n\in\mathbb N}\subset\calG.
\end{equation}
\item[(iii)] Restricted to $N_\calG$,
%
\begin{equation}\label{eqNG}
\sumn\widehat\phi_{u_n}f(s) < \infty,\qquad\mu\mbox{-a.e.}\mbox{ for some } \{
\phi
_{u_n}\}_{n\in\mathbb N}\subset\calG.
\end{equation}
\end{longlist}
\end{Thm}

The decomposition in (\ref{ePos-Null}) is unique (mod $\mu$). It is
referred to as the \textit{positive-null}
decomposition w.r.t. $\calG$. The sets $P_\calG$ and $N_\calG$ are
referred to as the \textit{positive} and \textit{null} parts of $S$ w.r.t.
$\calG$,
respectively. If $\mu(N_\calG) = 0$ [$\mu(P_\calG) = 0$, resp.], then
$\calG$ is said to be a \textit{positive} (\textit{null}, resp.)
$G$-action.

The next result provides an equivalent characterization of (\ref
{ePos-Null}) based on
the notion of a weakly wandering set. Recall that a measurable set
$W\subset S$ is \textit{weakly wandering}, w.r.t. $\calG$, if there
exists $\{\phi_{t_n}\}_{n\in\mathbb N}\subset\calG$ such that $\mu
(\phi
_{t_n}\inv(W)\cap\phi_{t_m}\inv(W)) = 0$
for all $n\neq m$.
%
\begin{Thm}\label{thmweaklyWandering}
Under the assumptions of Theorem~\ref{thmdecompPN}, we have the following:
\begin{longlist}[(ii)]
\item[(i)] The positive part $P_\calG$ has no weakly wandering set of
positive measure.
\item[(ii)] The null part $N_\calG$ is a union of weakly wandering sets
w.r.t. $\calG$.
\end{longlist}
\end{Thm}
%
\begin{Rem}\label{rem1} In the one-dimensional case, Krengel
\cite{krengel67classification} (for $G = \mathbb Z$) and
Samorodnitsky~\cite{samorodnitsky05null} (for $G = \mathbb R$)
establish alternative characterizations of
the decomposition (\ref{ePos-Null}). These results involve certain
integral tests, which we were unable to
extend to multiple dimensions. Takahashi's characterizations, employed
in Theorem~\ref{thmdecompPN}, are
valid for all dimensions.
\end{Rem}
%

\subsection{Multiparameter ergodic theorems}\label{secmultiparameter}

In the rest of the paper we focus on $\mathbb T^d$-actions, where
$\mathbb T$ stands for either $\mathbb Z$ or $\mathbb R$. We equip
$\mathbb T^d$ with the measure
\mbox{$\lambda\equiv\lambda_{\mathbb T^d}$}, which is either the counting
(if $\mathbb T=\mathbb Z$) or the Lebesgue (if $\mathbb T=\mathbb R$) measure.
In the sequel we establish multiparameter versions of the
\textit{stochastic ergodic theorem} and \textit{Birkhoff theorem} for the
case of $\mathbb T^d$-actions. They are extensions of the well-known
results in the one-dimensional case. The proofs follow from the works
of Krengel and Tempel'man (see, e.g.,~\cite{krengel85ergodic}).

Introduce the \textit{average functional} $A_T$, defined for all locally
integrable $h\dvtx\mathbb T^d\to\mathbb R$:
\[
A_Th \equiv A_{\mathbb T^d,T} h \defe\frac1{\ctt}\int_{\btt
}h(t)\lambda
(\d t)
\]
with $\btt\equiv B_{\mathbb T^d}(T) \defe(-T,T]^d\cap{\mathbb T}^d $
and $\ctt\equiv C_{\mathbb T^d}(T) \defe(2T)^d$.\vspace*{1pt}

Consider now a collection of functions $\{f_t\}_{t\in\mathbb
T^d}\subset\lone\smu$
such that
$(t,s)\mapsto f(t,s) \equiv f_t(s)$ is jointly measurable when
${\mathbb T}\equiv{\mathbb R}$. Then,
one can define the \textit{average operator}:
%
\begin{equation}\label{defavgop}
(A_Tf)(s) \defe\frac1{\ctt}\int_{ \btt} f_t(s)\lambda(\d t).
\end{equation}
Let \mbox{$\|\cdot\|$} denote the $L^1$ norm. If $t\mapsto\|
f_t\|$ is locally
integrable (i.e., integrable on finite intervals), then Fubini's
theorem implies that
$A_T f \in\lone\smu$, for all $T>0$.
Recall also that a sequence of measurable functions $\{f_n\}_{n\in
\mathbb N}\subset\la
\smu$ converges \textit{stochastically} (or
locally in measure) to $g\in\la\smu$, in short, $f_n\stackrel{\mu
}{\to} g$, as
$n\to
\infty$, if
%
\begin{eqnarray}\label{eqsfn}
&&\limn\mu\bigl(\{ s\dvtx |f_n(s)-g(s)|>\varepsilon\} \cap B\bigr) =
0\nonumber\\[-8pt]\\[-8pt]
&&\eqntext{\mbox{for all }
\varepsilon>0, B\in{\calB}\mbox{ with } \mu(B)<\infty.}
\end{eqnarray}

\begin{Rem}
By Theorem A.1 in~\cite{kolodynski03group}, there exists a strictly
positive measurable function $(t,s)\mapsto w(t,s)$, such that
for all $t\in{\mathbb T}^d$, $w(t,s) = \d(\mu\circ{\phi_t})/\d\mu(s)$
for $\mu$-almost all $s$, and for all $t,h \in\mathbb{T}^d$ and for
all $s \in S$,
%
\begin{equation}\label{ewt+h}
w(t+h,s) = w(h,s) w(t,\phi_h(s)).
\end{equation}
From now on, we shall use $w(t,s)$ as the version of the Radon--Nikodym
derivative $\d(\mu\circ\phi_t)/\d\mu(s)$.
\end{Rem}
%
\begin{Thm}[(Multiparameter stochastic ergodic theorem for nonsingular
actions)]
\label{thmMSET}
Let $\{\phi_t\}_{t\in\mathbb T^d}$ be a nonsingular $\mathbb
T^d$-action on the measure space $(S,\mu)$. Let $f_0\in\lone\smu$ and
define $f(t,s)\equiv(\what\phi_{-t} f_0)(s):=w(t,s)f_0\circ\phi_t(s)$.
Then, there exists $\widetilde f\in\lone\smu$, such that
%
\begin{equation}\label{eqassumption}
A_Tf \equiv\frac{1}{C(T)} \int_{B(T)} f(t,\cdot) \lambda(\d t)
\stackrel{\mu}{\to}
{\widetilde f}\qquad\mbox{as } T\to\infty.
\end{equation}
Moreover, $\widetilde f$ is invariant w.r.t. $\widehat\calG$, that is,
$\widehat\phi_t\widetilde f= \widetilde f$
for all $t\in\mathbb T^d$.\vadjust{\goodbreak}
\end{Thm}
\begin{pf}
\textit{Suppose first that ${\mathbb T} = {\mathbb Z}$.} The existence
of $\widetilde f$ follows from Krengel's stochastic ergodic theorem
(Theorem 6.3.10 in~\cite{krengel85ergodic}). To\vspace*{1pt} see that $\widetilde f$
is $L^1$-integrable, pick a subsequence $T_n$ such that
$A_{T_n}f\to\widetilde f,\mu$-a.e., as $n\to \infty$. By Fatou's lemma,
$\|\widetilde f\| = \|{\limn A_{T_n}f}\| \leq\liminf_{n\to\infty
}\|A_{T_n}f\|\le\|f_0\|<\infty$, which implies $\widetilde
f\in\lone\smu$. Here we used the fact that $\int_{S} |A_T f| \,\d\mu\le
A_T \int_{S} | \what\phi_{-t} f_0| \,\d\mu= A_T \|f_0\| =
\|f_0\|$.\vspace*{2pt}

We now prove that $\widetilde f$ is invariant w.r.t. $\widehat\calG$.
Fix $\tau\in{\mathbb T}^d$ and let
$T_n\to\infty$ be such that $g_n:=A_{T_n} f \to\wtilde f$, $\mu$-a.e.,
as $n\to\infty$. Then, since $\phi_{\tau}$
is nonsingular,
%
\begin{eqnarray}\label{ewhat-phi-f-tilde}
&&(\what\phi_{-\tau} g_n)(s) \equiv\frac{\d(\mu\circ\phi_{\tau
})}{\d\mu}(s)
g_n\circ\phi_\tau(s)\nonumber\\[-8pt]\\[-8pt]
&&\quad\longrightarrow\quad
\frac{\d(\mu\circ\phi_{\tau})}{\d\mu}(s) \wtilde f\circ\phi
_\tau(s)
\equiv
(\what\phi_{-\tau} \wtilde f)(s),\qquad \mu\mbox{-a.e.}
\nonumber
\end{eqnarray}
as $n\to\infty$. On the other hand, since $f(t,\phi_\tau(s)) =
w(t,\phi
_\tau(s))f_0\circ\phi_{t+\tau}(s)$,
we obtain by (\ref{ewt+h}) and Fubini's theorem that
\begin{eqnarray*}
(\what\phi_{-\tau} g_n)(s) &=& \frac{1}{C(T_n)} \int_{B(T_n)} w(\tau
+t,s) f_0(\phi_{\tau+t}(s)) \lambda(\d t)
\\
&=& \frac{1}{C(T_n)} \int_{B(T_n)+\tau} f(t,s) \lambda(\d t),
\qquad\mu\mbox{-a.e.}
\end{eqnarray*}
Therefore, by performing cancelations and applying Fubini's theorem, we~get
\[
\| \what\phi_{-\tau} g_n - g_n\| \le\frac{\lambda((B(T_n)+\tau
)\Delta
B(T_n))}{C(T_n)} \|f_0\|,
\]
where $D\Delta E = (D\setminus E)\cup(E\setminus D)$ is the symmetric
difference of sets. The last term vanishes,
as $n\to\infty$, since $\tau\in{\mathbb Z}^d$ is fixed. This implies
that $\what\phi_{-\tau} g_n \stackrel{\mu}{\to} \wtilde f$, as
$n\to\infty$, which, in view of (\ref{ewhat-phi-f-tilde}), yields
$\what\phi_{-\tau} \wtilde f = \wtilde f, \mu$-a.e.
This, since $\tau\in{\mathbb Z}^d$ was arbitrary, establishes the
desired invariance of the limit $\wtilde f$.

\textit{Suppose now that ${\mathbb T} = {\mathbb R}$.} Since we will use
the result proved for $\mathbb T = \mathbb Z$, we explicitly
write $A_{\mathbb Z^d,T}$ and $A_{\mathbb R^d,T}$ to distinguish
between the discrete and integral average operators, respectively.
In view of part (i), for all $\delta>0$, we have
%
\begin{eqnarray}\label{eqAndelta}
A_{\mathbb R^d,n\delta}
f_0 &\equiv& \frac1{(2n\delta)^d}\int
_{(-n\delta
,n\delta]^d} \widehat\phi_{-\tau} f\,\d\tau
\nonumber\\[-8pt]\\[-8pt]
&=& \frac{1}{(2n)^d} \sum_{t\in(-n,n]^d\cap\mathbb Z^d} \widehat\phi
_{-\delta t} g^{(\delta)} \equiv
A_{\mathbb Z^d,n} g^{(\delta)},
\nonumber
\end{eqnarray}
where
\[
g^{(\delta)}(s) \defe\frac{1}{\delta^d}\int_{(-\delta,0]^d}
(\what\phi
_{-\tau} f_0)(s)\,\d\tau\in\lone\smu.\vadjust{\goodbreak}
\]
As already shown for the case $\mathbb T = \mathbb Z$, the right-hand
side of (\ref{eqAndelta}) converges stochastically, as $n\to\infty$,
to $\widetilde g^{(\delta)}\in\lone\smu$, where $\wtilde g^{(\delta)}$
is $\what\phi_{-\delta t}$-invariant, for all $t\in{\mathbb Z}^d$.
Write $T_\delta= \lfloor T/\delta\rfloor\delta$. Since for all
$\delta
>0$, the volume of $(-T,T]^d\setminus(-T_\delta,T_\delta]^d$ is
$o(C(T))$ as $T\to\infty$, it follows that
$\|A_{\mathbb R^d,T}f-A_{\mathbb R^d,T_\delta}f\|\to0$ as $T\to
\infty$. Therefore, we have that
\[
A_{\mathbb R^d,T}f\stackrel{\mu}{\to}{\widetilde g}^{(\delta
)}\qquad\mbox{as } T\to
\infty,
\]
which shows, in particular, that $\widetilde g^{(\delta)}= \wtilde g
\in\lone\smu$ must be independent of $\delta>0$.
Since $\widetilde g$ is invariant w.r.t. $\widehat\phi_{\delta t}$ for
all $\delta>0$ and $t\in{\mathbb Z}^d$, it follows that
$\widetilde g$ is $\what{\cal G}$-invariant.
\end{pf}
%
\begin{Thm}[(Multiparameter Birkhoff theorem)]\label{thmMBT}
Assume the conditions of Theorem~\ref{thmMSET} hold. Suppose,
moreover, that the action $\{\phi_t\}_{t\in\mathbb T^d}$ is measure preserving
on $(S,\mu)$, and that $\mu$ is a probability measure. Then,
\[
A_Tf\to\widetilde f\defe\esp_\mu(f|\calI)
\qquad\mbox{almost surely
and in
}\lone,
\]
where $\calI$ is the $\sigma$-algebra of all $\calG$-invariant
measurable sets.
\end{Thm}
\begin{pf}
Suppose first that $\mathbb T = \mathbb Z$. The almost sure convergence
and the structure of the limit $\widetilde f$ follow
from Tempel'man's theorem (Theorem 6.2.8 in~\cite{krengel85ergodic},
page 205). The $\lone$-convergence is clear when $f_0$ is bounded.
Suppose now that
$f_0\in\lone\smu$. Consider the sequence $A_T f, T\in\bbN$. For all
$\varepsilon>0$ there exists a bounded $f_0^{(\varepsilon)}\in L^\infty
\smu$
such that
$\|f_0-f_0^{(\varepsilon)}\|<\varepsilon/3$. Then, by the triangle
inequality and the fact that $A_T$ is a linear contraction, we~get
\begin{eqnarray*}
\|A_{T_1} f - A_{T_2}f\| &\le& \bigl\|A_{T_1} f^{(\varepsilon)} - A_{T_2}
f^{(\varepsilon)}\bigr\| + 2\bigl\| f_0-f_0^{(\varepsilon)}\bigr\| \\
&\le& \bigl\|A_{T_1} f^{(\varepsilon)} - A_{T_2} f^{(\varepsilon)}\bigr\|
+ 2\varepsilon/3 < \varepsilon
\end{eqnarray*}
for all sufficiently large $T_1$ and $T_2$. This is because $A_T
f^{(\varepsilon)}$ converges in $L^1$. We have thus shown
that $A_Tf, T\in\mathbb N$, is a Cauchy sequence in the Banach space
$L^1(S,\mu)$ and, hence, it has a limit, which is necessarily
$\wtilde f$.

Let now $\mathbb T = \mathbb R$. First, by a discretization argument as
in the proof of Theorem~\ref{thmMSET}, we can show
$A_Tf \to\widetilde f$ almost surely, for all $f_0\in L^1(S,\mu)$. The
$L^1$-convergence can be established as in the proof in the discrete case.
\end{pf}

\section{Stationary sum-stable random fields}\label{secdecomp}

We focus on S$\alpha$S ($0<\alpha<2$) random fields $\vv X=\{X_t\}
_{t\in\mathbb T^d}$, with a \textit{spectral representation}:
%
\begin{equation}\label{repintegralRep}
\{X_t\}_{t\in\mathbb T^d} \eqd\biggl\{\int_S f_t(s) M_\alpha(\d
s)\biggr\}_{t\in \mathbb T^d}.
\end{equation}
Here $\{f_t\}_{t\in\mathbb T^d}\subset\la(S,\mu)$, and the integral
is with respect to
an independently scattered
\textit{S$\alpha$S random measure} $M_\alpha$ on $S$ with control measure
$\mu$
(see Chapters~3 and 13 in~\cite{samorodnitsky94stable} for more details).
Without loss of generality, we shall also assume that $\{f_t\}_{t\in
\mathbb T^d}$ has
\textit{full support} in $\la\smu$. Namely,
there is no $B\in\calB$ with \mbox{$\mu(B)>0$}, such that $\int
_B|f_t(s)|^\alpha\mu(\d s) = 0$, for all $t\in\mathbb T^d$.

All measurable S$\alpha$S random fields $\vv X$ have a spectral
representation (\ref{repintegralRep}), where $(S,\mu)$ can be chosen
to be a {standard Lebesgue space}
and the functions $(t,s)\mapsto f_t(s)$ to be jointly measurable (see,
e.g., Proposition 11.1.1 and Theorem 13.2.1 in~\cite{samorodnitsky94stable}).

It is known from Rosi\'nski~\cite{rosinski95structure,rosinski00decomposition}
that when $\vv X$ is stationary,
there exists a \textit{minimal} spectral representation (\ref
{repintegralRep}) with
%
\begin{equation}\label{rep2}
f_t(s) = c_t(s)\biggl(\frac{\d(\mu\circ\phi_t)}{\d\mu}(s)
\biggr)^{1/\alpha
}f_0\circ
\phi_t(s),\qquad t\in\mathbb T^d,
\end{equation}
where $f_0\in\la\smu$, $\{\phi_t\}_{t\in\mathbb T^d}$ is a
nonsingular $\mathbb
T^d$-action on $\sbm$, and $\{c_t\}_{t\in\mathbb T^d}$ is a cocycle
for $\{\phi_t\}_{t\in\mathbb T^d}$ taking values in $\{-1,1\}$.
Namely, $(t,s)\mapsto
c_t(s) \in\{-1,1\}$ is a measurable map,
such that for all $u,v\in\mathbb T^d$, $c_{u+v}(s) = c_v(s)c_u
(\phi
_v(s)),\mu$-a.e. $s\in S$.
The representation (\ref{repintegralRep}) is minimal, if the ratio
$\sigma$-algebra $\sigma(f_t/f_\tau\dvtx t,\tau\in\mathbb T^d)$ is
equivalent to $\calB$ (see Definition 2.1 in~\cite{rosinski95structure}).
The minimality is an indispensable tool to study the spectral
representations, although it is hard to check in practice. For more
equivalent conditions and insights, see Rosi\'nski
\cite{rosinski06minimal} and Pipiras~\cite{pipiras07nonminimal}.

We say that a random field $\{X_t\}_{t\in\mathbb T^d}$ with the minimal
representation (\ref{repintegralRep}) and (\ref{rep2})
is \textit{generated by the $\mathbb T^d$-action} $\{\phi_t\}_{t\in
\mathbb T^d}$ \textit{and
the cocycle $\{c_t\}_{t\in\mathbb T^d}$}. In this case, we also say $\{
X_t\}_{t\in\mathbb T^d}$ has an action
representation $(f_0,\calG\equiv\{\phi_t\}_{t\in\mathbb T^d}, \{
c_t\}_{t\in\mathbb T^d})$.

It turns out, moreover, the action $\{\phi_t\}_{t\in\mathbb T^d}$ is
determined by the
distribution of $\{X_t\}_{t\in\mathbb T^d}$, up to the equivalence
relationship of
$\mathbb T^d$-actions (see Theorem 3.6 in~\cite{rosinski95structure}).
Thus, structural results for the $\mathbb T^d$-actions
imply important structural results for the corresponding S$\alpha$S
random fields. In particular, by using
Theorem~\ref{thmdecompPN}, we obtain the following result:
%
\begin{Thm}\label{thm1}
Let $\{X_t\}_{t\in\mathbb T^d}$ be a measurable stationary S$\alpha$S
random field with
spectral representation
(\ref{repintegralRep}). We suppose that $(S,{\mathcal B},\mu)$ is a
standard Lebesgue space and the spectral
representation $\{f_t(s)\}_{t\in\mathbb T^d}$ is measurable. Assume,
in addition, that
%
\begin{equation}\label{eqgs}\quad
g(s) \defe\int_{T_0}a_\tau|f_{\tau}(s)|^\alpha\lambda(\d\tau)
\mbox{
is } L^1\mbox{-integrable}\quad\mbox{and}\quad\supp(g) = S
\end{equation}
for some $T_0\in{\mathcal B_{\mathbb T^d}}$ and $a_\tau>0, \forall
\tau\in T_0$. Then:
\begin{longlist}[(ii)]
\item[(i)] $\{X_t\}_{t\in\mathbb T^d}$ is generated by a positive
$\mathbb T^d$-action
if and only if
%
\begin{equation}\label{eqPG1}\quad
\sumn\int_{T_0}a_\tau|f_{\tau+ t_n}(s)|^\alpha\lambda(\d\tau) =
\infty
,\qquad\mu\mbox{-a.e.}\mbox{ for all }
\{t_n\}_{n\in\mathbb N}\subset\mathbb T^d.
\end{equation}
\item[(ii)] $\{X_t\}_{t\in\mathbb Z^d}$ is generated by a null
$\mathbb T^d$-action if
and only if
%
\begin{equation}\label{eqNG1}\qquad
\sumn\int_{T_0}a_\tau|f_{\tau+ t_n}(s)|^\alpha\lambda(\d\tau) <
\infty
,\qquad\mu\mbox{-a.e.}\mbox{ for some }
\{t_n\}_{n\in\mathbb N}\subset\mathbb T^d.
\end{equation}
\end{longlist}
In particular, the classes of stationary S$\alpha$S random fields
generated by positive and null
$\mathbb T^d$-actions are disjoint.
\end{Thm}
%
\begin{Rem}
One can always choose $\{a_\tau\}_{\tau\in T_0}$ such that (\ref
{eqgs}) holds, if the spectral functions $\{f_t\}_{t\in\mathbb T^d}$ have
full support in $\la\smu$.
\end{Rem}
\begin{pf*}{Proof of Theorem~\ref{thm1}}
Suppose first that $\{f_t\}_{t\in\mathbb T^d}$ is minimal and, hence,
it has the
form (\ref{rep2}).
Observe that, for all $t,\tau\in\mathbb T^d$, we have
\[
|f_{\tau+t}(s)|^\alpha= \frac{\d(\mu\circ\phi_{t})}{\d\mu}(s)
\,\frac{\d(\mu\circ\phi_{\tau})}{\d\mu}\circ\phi_{t}(s)
|f_0\circ\phi_{\tau}\circ\phi_t(s)|^\alpha, \qquad\mu\mbox{-a.e.}
\]
Since both the left-hand side and the right-hand side are measurable in $(\tau,s)$, by
Fubini's theorem,
\begin{eqnarray*}
\int_{T_0} a_\tau|f_{\tau+t}(s)|^\alpha\lambda(\d\tau)
&=& \frac{\d(\mu\circ\phi_{t})}{\d\mu}(s) \int_{T_0}
a_\tau|f_{\tau}\circ\phi_t(s)|^\alpha\lambda(\d\tau)\\
&=& (\what
\phi
_{-t} g)(s), \qquad\mu\mbox{-a.e.},
\end{eqnarray*}
where the last relation follows from (\ref{eqdual}).
Therefore,
\[
\sumn\int_{T_0}a_\tau|f_{\tau+t_n}(s)|^\alpha\lambda(\d\tau) =
\sumn\what\phi_{-t_n} g, \qquad\mu\mbox{-a.e. } \forall\{t_n\}_{n\in\mathbb N}
\subset{\mathbb T}^d.
\]
Hence, Theorem~\ref{thmdecompPN}(ii) and (iii), applied
to the strictly positive function $g\in L^1(S,\mu)$,
implies the statements of parts (i) and (ii), respectively.

Using Remark 2.5 in~\cite{rosinski95structure} and a standard Fubini
argument, it can be shown that a test function (\ref{eqgs}) in the
general case corresponds to one in the situation when the integral
representation $\{f_t\}_{t\in\mathbb T^d}$ of the field is of the
form (\ref{rep2}). Therefore, an argument parallel to the proof of
Corollary 4.2 in~\cite{rosinski95structure} shows that the tests
described in this theorem can be applied to any full support integral
representation, not necessarily minimal or of the form (\ref{rep2}).
This completes the proof.
\end{pf*}

The above characterization motivates the following decomposition of an
arbitrary measurable stationary S$\alpha$S
random field $\vv X=\{X_t\}_{t\in\mathbb T^d}$. Without loss of
generality, let $\vv X$
have a representation $(f_0,\calG\equiv\{\phi_t\}_{t\in\mathbb T^d},
\{c_t\}_{t\in\mathbb T^d})$ as in (\ref{repintegralRep}) and~(\ref
{rep2}). Then, by Lemma~\ref{lemmaximal},
$S = P_\calG\cup N_\calG$ and one can write
%
\begin{equation}\label{decompPN}
\{X_t\}_{t\in\mathbb T^d} \eqd\{X^P_t+X^N_t\}_{t\in\mathbb T^d}
\end{equation}
with
\[
X_t^P = \int_{P_\calG}f_t(s)M_\alpha(\d s) \quad\mbox{and}\quad
X_t^N = \int
_{N_\calG}f_t(s)M_\alpha(\d s)\qquad\mbox{for all }t\in\mathbb T^d.\vadjust{\goodbreak}
\]

\begin{Coro}\label{corodecompPN} \textup{(i)} The decomposition (\ref
{decompPN}) is unique in law. That is, if there is another
representation $(f_0^{(2)},\calG^{(2)}\equiv\{\phi^{(2)}_t\}
_{t\in\mathbb T^d},\{c^{(2)}_t\}_{t\in\mathbb T^d})$
satisfying (\ref{repintegralRep}) and (\ref{rep2}), then
\[
\{X^P_t\}\stackrel{\mathit{d}}{=}\biggl\{\int_{P_{\calG^{(2)}}}f^{(2)}_t\,\d M_\alpha
\biggr\}
\quad\mbox{and}\quad\{X^N_t\}\stackrel{\mathit{d}}{=}
\biggl\{\int_{N_{\calG^{(2)}}}f^{(2)}_t\,\d
M_\alpha
\biggr\}.
\]

\textup{(ii)} The components $\bX^P=\{X^P_t\}_{t\in\mathbb T^d}$ and
$\bX^N=\{X^N_t\}_{t\in\mathbb T^d}$ are independent, $\bX^P$~is
generated by a positive $\mathbb
T^d$-action and
$\bX^N$ is generated by a null $\mathbb T^d$-action.
\end{Coro}
\begin{pf} Proof of (ii) is trivial. To prove (i), observe that by
Remark 2.5 in~\cite{rosinski95structure}, there exist measurable functions
$\Phi\dvtx S_2\to S$ and $h\dvtx S_2\to\mathbb R\setminus\{0\}$ such that for all
$t\in\mathbb T^d$,
%
\begin{equation}\label{eqf2}
f^{(2)}_t(s) = h(s) f_t\circ\Phi(s), \qquad\mu_2\mbox{-almost all }
s\in
S_2
\end{equation}
and\vspace*{1pt} $\d\mu= (|h|^\alpha\,\d\mu_2)\circ\Phi\inv$. Using
(\ref {eqf2}) and an argument parallel to the proof of (2.18) in
\cite{samorodnitsky05null}, it can be shown that $P_{\calG^{(2)}}=\Phi
^{-1}(P_{\calG})$ and $N_{\calG^{(2)}}=\Phi^{-1} (N_{\calG})$ modulo
$\mu_2$, from which the distributional equality in (i) follows as in
the proof of Theorem 4.3 in~\cite{rosinski95structure}.
\end{pf}



\section{\texorpdfstring{Ergodic properties of stationary S$\alpha$S fields}
{Ergodic properties of stationary S alpha S fields}}\label{secergo}
Let $(\Omega,\filF,\proba)$ be a probability space, and $\{\theta
_t\}_{t\in\mathbb T^d}$
a measure-preserving $\mathbb T^d$-action on
$(\Omega,\filF,\proba)$. Consider the random field $X_t(\omega) =
X_0\circ\theta_t(\omega)$, $t\in\mathbb T^d$. The
random field $\{X_t\}_{t\in\mathbb T^d}$ defined in this way is
stationary and,
conversely, any stationary measurable
random field can be expressed in this form.

\textit{We start by introducing some notation.} For $t\in\mathbb T^d$,
let $\|t\|$ denote its sup norm.
We consider the class $\calT$ of all sequences that converge to infinity:
\[
\calT\defe\Bigl\{\{t_n\}_{n\in\mathbb N}\subset\mathbb T^d\dvtx\limn
\|t_n\| = \infty\Bigr\}.
\]
Recall that a set $E\subset\mathbb T^d$ is said to have \textit{density
zero} in $\mathbb T^d$ if
%
\begin{equation}\label{eqdensity0}
\lim_{T\to\infty}\frac1{\ctt}\int_{\btt}\ind_E(t)\lambda(\d t)
= 0.
\end{equation}
A set $D\subset\mathbb T^d$ is said to have \textit{density one} in
$\mathbb T^d$ if $\mathbb T^d\setminus D$ has density zero in $\mathbb
T^d$. The class of all sequences on $D$ that converge to infinity will
be denoted by
\[
\calT_D \defe\Bigl\{\{t_n\}_{n\in\mathbb N}\dvtx t_n\in\mathbb T^d\cap
D,\limn\|t_n\| = \infty\Bigr\}.
\]

\textit{Now we recall some basic definitions.} Write $\sigma_{\vv X}
\defe
\sigma(\{X_t\dvtx t\in\mathbb T^d\})$ for the $\sigma$-algebra generated by
the field $\{X_t\}_{t\in\mathbb T^d}$.
We say $\{X_t\}_{t\in\mathbb T^d}$ is:
\begin{longlist}[(iii)]
\item[(i)] \textit{ergodic}, if
%
\begin{equation}\label{eqergodic}
\limT\frac1{\ctt}\int_{\btt}\proba\bigl(A\cap\theta_t(B)\bigr)\lambda(\d
t) =
\proba(A)\proba(B) \qquad\mbox{for all } A,B\in\sigma_{\vv X}.\hspace*{-32pt}
\end{equation}

\item[(ii)] \textit{weakly mixing}, if there exists a density one set $D$
such that
%
\begin{equation}\label{eqweakMixing}
\limn\proba\bigl(A\cap\theta_{t_n}(B)\bigr) = \proba(A)\proba(B)
\qquad\mbox{for all
} A,B\in\sigma_{\vv X},\{t_n\}_{n\in\mathbb N}\in\calT_D.\hspace*{-32pt}
\end{equation}
\item[(iii)] \textit{mixing}, if
%
\begin{equation}\label{eqmixing}
\limn\proba\bigl(A\cap\theta_{t_n}(B)\bigr) = \proba(A)\proba(B)
\qquad\mbox{for all
} A,B\in\sigma_{\vv X},\{t_n\}_{n\in\mathbb N}\in\calT.\hspace*{-32pt}
\end{equation}
\end{longlist}
In general, we always have that
\[
\mbox{mixing} \Rightarrow\mbox{weakly mixing} \Rightarrow\mbox{ergodicity.}
\]
For stationary S$\alpha$S random fields, however, we have the
following result.
%
\begin{Thm}\label{thmsumErgodicity}
Let $\{X_t\}_{t\in\mathbb T^d}$ denote a measurable S$\alpha$S random
field with spectral
representation (\ref{repintegralRep}) and $\alpha\in(0,2)$.
The following are equivalent:

\begin{longlist}[(iii)]
\item[(i)] $\{X_t\}_{t\in\mathbb T^d}$ is ergodic.

\item[(ii)] $\{X_t\}_{t\in\mathbb T^d}$ is weakly mixing.

\item[(iii)] $\lim_{T\to\infty}C(T)\inv\int_{B(T)}\exp(2\|f_0\|
_\alpha^\alpha- \|f_0-f_t\|_\alpha^\alpha)\lambda(\d t) = 1$.

\item[(iv)] The $\mathbb T^d$-action $\{\phi_t\}_{t\in\mathbb
T^d}$ has no nontrivial
positive component.
\end{longlist}
\end{Thm}
\begin{pf}
Using Theorem~\ref{thmMBT} and proceeding as in Theorems 2 and 3
in~\cite{podgorski92note}, one can show the equivalence of (i), (ii)
and (iii).

To prove the equivalence of (ii) and (iv), we need the following
result, which is an extension of Theorem 2.7 in~\cite{gross94some}. The
proof is given in the \hyperref[app]{Appendix}. We also fill a gap in the results of
Gross~\cite{gross94some} (see Remark~\ref{remarkmistakegross}).
%
\begin{Prop}\label{propgross}
Assume $\alpha\in(0,2)$ and $\{X_t\}_{t\in\mathbb T^d}$ is a
stationary S$\alpha$S
random field with spectral representation $\{f_t\}_{t\in\mathbb
T^d}\subset\la\sbm$.
Then, the process $\{X_t\}_{t\in\mathbb T^d}$ is weakly mixing if and
only if there
exists a density one set $D\subset\mathbb T^d$, such that
%
\begin{eqnarray}\label{eqgross94}
\limn\mu\{s\dvtx|f_0(s)|^\alpha\in K, |f_{t^*_n}(s)|^\alpha>
\varepsilon\} = 0\nonumber\\[-8pt]\\[-8pt]
&&\eqntext{\mbox{for all compact } K\subset\mathbb R\setminus\{0\}, \varepsilon
>0 \mbox{ and } \{t^*_n\}_{n\in\mathbb N}\in\calT_D.}
\end{eqnarray}
\end{Prop}

Now we prove the equivalence of (ii) and (iv) by following closely the
proof of Theorem 3.1 in~\cite{samorodnitsky05null}. The proof of (ii)
implying (iv) remains the same. To show that (iv) implies (ii),
however, we treat the discrete and the continuous parameter scenarios
together by virtue of Theorem~\ref{thmMSET}, which unifies the two
cases (which were treated differently in~\cite{samorodnitsky05null}).
More specifically, in view of (\ref{eqgross94}) and a multivariate
extension of Lemma 6.2 in~\cite{petersen83ergodic}, page 65, it is
enough to show that for all $\varepsilon>0$ and compact sets $K\subset
\mathbb R\setminus\{0\}$,
%
\begin{equation}\label{eqfrac12T}
\limT A_T\mu\bigl\{s\dvtx|f_0(s)|^\alpha\in K,\bigl|f_{(\cdot)}(s)\bigr|^\alpha
>\varepsilon\bigr\} = 0,
\end{equation}
where $A_T$ is the average operator defined by (\ref{defavgop}).
Following verbatim the argument in the proof of $(3.1)$ in \cite
{samorodnitsky05null}, we obtain (\ref{eqfrac12T}) for both discrete
and continuous parameter cases with the help of Theorem~\ref{thmMSET}.
\end{pf}
%
\begin{Rem}
From the structural results~\cite{roy08stationary,roy10nonsingular} and
Theorem~\ref{thmsumErgodicity} above, we obtain a unique in law
decomposition of $\mathbf{X}$ into three independent stable processes
in parallel to the one-dimensional case~\cite{samorodnitsky05null},
that is,
\[
\mathbf{X}=\mathbf{X}^{(1)}+\mathbf{X}^{(2)}+\mathbf{X}^{(3)},
\]
where $\mathbf{X}^{(1)}$ is a mixed moving average in the sense of
\cite{surgailis93stable}, $\mathbf{X}^{(2)}$ is weakly mixing with no mixed
moving average component and $\mathbf{X}^{(3)}$ has no weakly mixing component.
\end{Rem}

\section{Max-stable stationary random fields} \label{secmax-stable}
In this section we discuss the structure and ergodic properties of
stationary max-stable random fields, indexed by
$\mathbb T^d$. For simplicity and without loss of generality, we will
focus on $\alpha$-Fr\'echet random fields. The random
field $\vv Y = \{Y_t\}_{t\in\mathbb T^d}$ is said to be
$\alpha$-Fr\'echet, if for all $a_j>0, \tau_j\in\mathbb T^d,
1\le j\le n$, the max-linear combinations $\xi:= \max_{1\le j\le n} a_j
Y_{\tau_j} \equiv\bigvee_{1\le j\le n} a_j Y_{\tau_j}$ have $\alpha
$-Fr\'echet distributions. Namely,
$\proba(\xi\leq x) = \exp\{-\sigma^\alpha x^{-\alpha}\}\mbox{ for all
} x\in(0,\infty)$,
where $\sigma>0$ is referred to as the \textit{scale coefficient} and
$\alpha>0$ is the tail index of $\xi$.
The $\alpha$-Fr\'echet random fields are max-stable. Conversely, all
max-stable random fields with $\alpha$-Fr\'echet
marginals are $\alpha$-Fr\'echet random fields.

The spectral representations for $\alpha$-Fr\'echet random fields have
been developed by de Haan~\cite{dehaan84spectral} and developed
by \cite
{stoev06extremal,wang10structure}.
Any measurable $\alpha$-Fr\'echet random field $\vv Y = \{Y_t\}_{t\in
\mathbb T^d}$ ($\alpha>0$) can be
represented as
%
\begin{equation}\label{repextremalRep}
\{Y_t\}_{t\in\mathbb T^d} \eqd\biggl\{\Eintt_S f_t(s) \mam(\d s)
\biggr\}_{t\in\mathbb T^d},
\end{equation}
where $\{f_t\}_{t\in\mathbb T^d}\subset\lap(S,\mu)\defe\{f\in\la
(S,\mu)\dvtx f\geq0\}$,
``\mbox{$\Einttt$}'' stands for the \textit{extremal
integral}, $\mam$ is an
\textit{independently scattered $\alpha$-Fr\'echet random sup-measure}
with control measure $\mu$ and $(S,\mu)$ can be chosen to be a standard
Lebesgue space (see~\cite{stoev06extremal,wang10structure}). The
functions $\{f_t\}_{t\in\mathbb T^d}$ in (\ref{repextremalRep}) are
called \textit{spectral
functions} of the $\alpha$-Fr\'echet random field.
If the representation in (\ref{repextremalRep}) is \textit{minimal}, as
in the sum-stable case, it then follows that
%
\begin{equation}\label{eft-max-stable}
f_t(s) = \biggl(\frac{\d(\mu\circ\phi_t)}{\d\mu}\biggr)^{1/\alpha
}f_0\circ
\phi
_t(s)\qquad\mbox{for all } t\in\mathbb T^d,
\end{equation}
where $\phi=\{\phi_t\}_{t\in\mathbb T^d}$ is a nonsingular group
action and $f_0\in\lap(S,\mu)$ (see, e.g.,~\cite
{wang10structure}, Theorems 3.1 and 3.2).
Thus, the $\alpha$-Fr\'echet random field $\vv Y$ is said to be
generated by the group action $\phi$ if
(\ref{repextremalRep}) is a minimal representation such that (\ref
{eft-max-stable}) holds. This allows us to extend
the available classification results in the sum-stable case to the
max-stable setting.
Note that compared to (\ref{rep2}), the cocycle $\{c_t\}_{t\in\mathbb
T^d}$ disappears,
as $\{f_t\}_{t\in\mathbb T^d}$ are nonnegative.
By a similar argument as in Theorem~\ref{thm1}, we obtain the
following result.
%
\begin{Thm}\label{thmmax1}
Suppose $\{Y_t\}_{t\in\mathbb T^d}$ is a measurable stationary $\alpha
$-Fr\'echet random
field with spectral representation $\{f_t\}_{t\in\mathbb T^d}$
as in (\ref{repextremalRep}). Let $T_0 \in{\cal B}_{\mathbb T^d}$ and
$\{a_\tau\}_{\tau\in T_0}$, $a_\tau>0$, be
such that (\ref{eqgs}) holds. Then:
\begin{longlist}[(ii)]
\item[(i)] $\{Y_t\}_{t\in\mathbb T^d}$ is generated by a positive
$\mathbb T^d$-action,
if and only if (\ref{eqPG1}) holds.
\item[(ii)] $\{Y_t\}_{t\in\mathbb T^d}$ is generated by a null
$\mathbb T^d$-action, if
and only if (\ref{eqNG1}) holds.
\end{longlist}
In particular, the classes of stationary $\alpha$-Fr\'echet random
fields generated by positive and null $\mathbb T^d$-actions are disjoint.
\end{Thm}

An intimate connection between the $\alpha$-Fr\'echet and S$\alpha$S
processes\break \mbox{($0<\alpha<2$)} was recently revealed through the notion of
\textit{association}, independently by Kabluchko \cite
{kabluchko09spectral} and Wang and Stoev~\cite{wang10association}.
By the association tool established in~\cite{wang10association}, the
decomposition results for $\alpha$-Fr\'echet random fields follow
immediately from the corresponding ones for S$\alpha$S random fields.
Indeed, for an $\alpha$-Fr\'echet random field $\{Y_t\}_{t\in\mathbb
T^d}$ with spectral
functions $\{f_t\}_{t\in\mathbb T^d}$, $\alpha\in(0,2)$, consider
the S$\alpha$S random
field $\{X_t\}_{t\in\mathbb T^d}$ with the same spectral functions.
Naturally, the
random fields $\{X_t\}_{t\in\mathbb T^d}$ and $\{Y_t\}_{t\in\mathbb
T^d}$ are said to be \textit{associated}, according to~\cite{wang10association}. Then, applying
Theorem 5.1 in~\cite{wang10association} to Corollary~\ref
{corodecompPN}, we obtain the following results on $\alpha$-Fr\'echet
random fields.
%
\begin{Coro}
Let $\{Y_t\}_{t\in\mathbb T^d}$ be a measurable stationary $\alpha
$-Fr\'echet random
field with representation in the form of (\ref{repextremalRep})
and (\ref{eft-max-stable}).
We have the unique-in-law decomposition $\{Y_t\}_{t\in\mathbb T^d}
\stackrel{\mathit{d}}{=}\{Y^P_t \vee
Y^N_t\}_{t\in\mathbb T^d}$,
with
\[
Y_t^P = \Einttt_{P_\calG}f_t(s)\mam(\d s)\quad\mbox{and}\quad
Y_t^N = \Einttt
_{N_\calG
}f_t(s)\mam(\d s)\qquad\!\!\mbox{for all } t\in\mathbb T^d
\]
with\vspace*{1pt} $\calG\equiv\{\phi_t\}_{t\in\mathbb T^d}$.
The two components are independent, $\{Y^P_t\}_{t\in\mathbb T^d}$ is
generated by
positive $\mathbb T^d$-action and $\{Y^N_t\}_{t\in\mathbb T^d}$ is
generated by null
$\mathbb T^d$-action.
\end{Coro}

The ergodic properties of stationary $\alpha$-Fr\'echet random fields
can be characterized in terms of the recurrence properties of the
nonsingular group actions, as in the sum-stable case.
The following theorem extends the known results in the one-dimensional
case (see \cite
{stoev08ergodicity,kabluchko09spectral,kabluchko10ergodic}). These
results, however, cannot be established by
the association method.
%
\begin{Thm}\label{thmmaxErgodicity1}
Let $\{Y_t\}_{t\in\mathbb T^d}$ denote a measurable
$\alpha$-Fr\'echet random field
with spectral
representation (\ref{repextremalRep}) and (\ref{eft-max-stable}).
The following are equivalent:

\begin{longlist}[(iii)]
\item[(i)]
$\{Y_t\}_{t\in\mathbb T^d}$ is ergodic.

\item[(ii)]
$\{Y_t\}_{t\in\mathbb T^d}$ is weakly mixing.

\item[(iii)]
$\limT\ctt\inv\int_{\btt}\|f_t\wedge f_0\|_\alpha^\alpha
\lambda(\d
t) = 0$.

\item[(iv)] The $\mathbb T^d$-action $\{\phi_t\}_{t\in\mathbb
T^d}$ has no nontrivial
positive component.
\end{longlist}
\end{Thm}
\begin{pf}
The equivalence of (i), (ii) and (iii) for $\mathbb R$-action is proved
by Kabluchko and Schlather~\cite{kabluchko10ergodic}, Theorem 1.2.
Their proof generalizes to $\mathbb T^d$-actions as well.
The equivalence of (i) and (iv) can be proved by extending the proof of
Theorem 8 in~\cite{kabluchko09spectral} to the multiparameter setting,
using Theorems~\ref{thmMSET} and~\ref{thmMBT} accordingly.
\end{pf}


\section{Examples}\label{secexample}
This section contains two examples of stable random fields and their
ergodic properties via the positive-null decomposition of the
underlying action. These examples show the usefulness of our results to
check whether or not a stationary S$\alpha$S (or max-stable) random
field is ergodic (or, equivalently, weakly mixing).

The first example is based on a self-similar S$\alpha$S processes
with stationary increments introduced by~\cite{cohen06random} as a
stochastic integral with respect to an S$\alpha$S random measure, with
the integrand being the local time process of a fractional Brownian
motion. We extend these processes by replacing the fractional Brownian
motion by a Brownian sheet. We can call it a \textit{Brownian sheet local
time fractional S$\alpha$S random field} following the terminology of
\cite{cohen06random}.

\begin{Example} Suppose $(\Omega^\prime, \mathcal
{F}^\prime
, P^\prime)$ is a probability space supporting a Brownian sheet $\{
B_u\}
_{u \in\rpd}$. By~\cite{ehm81sample}, $\{B_u\}$ has a jointly
continuous local time field $\{l(x,u)\dvtx x \in\mathbb{R}, u \in\rpd
\}
$ defined on the same probability space. We will define an S$\alpha$S
random field based on this local time field, which inherits the
stationary increments property from $\{B_u\}_{u \in\rpd}$. Let
$M_\alpha$ be an S$\alpha$S random measure on $\Omega^\prime\times
\mathbb{R}$ with control measure $P^\prime\times \mathrm{Leb}$ living on
another probability space $(\Omega, \mathcal{F}, P)$. Following
verbatim the calculations of~\cite{cohen06random}, we have that
\[
Z_u = \int_{\Omega^\prime\times\mathbb{R}} l(x,u)(\omega^\prime
)M_\alpha(\d\omega^\prime,\d x),\qquad u \in\rpd,
\]
is a well-defined S$\alpha$S random field, which has stationary
increments over $d$-dimensional rectangles.

We now concentrate on the increments of $\{Z_u\}$ taken
over $d$-dimensional rectangles. For any $t \in\zpd$, define
%
\begin{equation}\label{defnX}
X_t=\Delta Z_t:=\sum_{i_1=0}^1\sum_{i_2=0}^1 \cdots\sum_{i_d=0}^1
(-1)^{i_1+i_2+\cdots+i_d+d}Z_{t+(i_1,i_2,\ldots,i_d)}.
\end{equation}
Clearly, $\{X_t\}_{t \in\zpd}$ is a stationary S$\alpha$S random
field, which can be extended (in law) to a stationary S$\alpha$S random
field $\mathbf{X}:=\{X_t\}_{t \in\Zd}$ by Kolmogorov's extension
theorem. We claim that $\mathbf{X}$ is generated by a null $\Zd
$-action. To prove this,\vspace*{1pt} define, for all $n \geq1$, $\tau
^{(n)}:=(n^{4/d}, n^{4/d}, \ldots, n^{4/d})$, and for all $n \geq1$
and $t \in\zpd$,
\[
T_{n,t}:=\{s\dvtx t_i+n^{4/d} \leq s_i \leq1+t_i+n^{4/d} \mbox{ for
all }i=1,2,\ldots,d\}.
\]
For each $t\in\zpd$, take a positive real number $a_t$ in such a way
that $\sum_{t \in\zpd} a_t=1$. Defining $\Delta l(x,t)$ in parallel to
(\ref{defnX}) and following the proof of $(4.7)$ in~\cite
{cohen06random}, we can establish that
\begin{eqnarray*}
&&\int_{\Omega^\prime}\int_{\mathbb{R}} e^{-x^2/2}\sum_{t \in\zpd
}\sum
_{n=1}^{\infty} a_t\Delta l\bigl(x,t+\tau^{(n)}\bigr)\,\d x\,\d
P^\prime\\
&&\qquad=\sum_{t \in\zpd}a_t\sum_{n=1}^{\infty} \int_{T_{n,t}} \frac{\d
s}{\sqrt{1+\prod_{i=1}^d s_i}} \leq\sum_{n=1}^{\infty} \frac
{1}{\sqrt
{1+n^4}} < \infty.
\end{eqnarray*}
This shows, in particular, that $\sum_{t \in\zpd}\sum_{n=1}^{\infty}
a_t\Delta l(x,t+\tau^{(n)})(\omega^\prime)<\infty$ for
$P^\prime\times \mathrm{Leb}$-almost all $(\omega^\prime,x) \in\Omega\times
\mathbb{R}$. Besides, it can be easily shown that $\sum_{t \in\zpd
}a_t\Delta l(x,\allowbreak t)(\omega^\prime)>0$ for $P^\prime\times
\mathrm{Leb}$-almost all $(\omega^\prime,x) \in\Omega\times\mathbb{R}$ (see,
e.g.,~\cite{tran76problem}). Hence, by Theorem~\ref{thm1}, it
follows that $\mathbf{X}$ is generated by a null action and hence is
weakly mixing.
\end{Example}

The next example is based on a class of mixing stationary S$\alpha$S process
considered in~\cite{rosinski96classes}. We look at a stationary
S$\alpha
$S random field generated by $d$ independent recurrent Markov chains,
at least one of which is null-recurrent. This is a class of stationary
S$\alpha$S random fields which are weakly mixing as a field but not
necessarily ergodic in every direction.
%
\begin{Example}
We start with $d$ irreducible aperiodic recurrent
Markov chains on $\mathbb{Z}$ with laws $P^{(1)}_i(\cdot),
P^{(2)}_i(\cdot), \ldots, P^{(d)}_i(\cdot)$, $i \in\mathbb{Z}$ and
transition probabilities $(p^{(1)}_{jk}),(p^{(2)}_{jk}), \ldots
,(p^{(d)}_{jk})$, respectively. For all $l=1,2,\ldots,d$, let
$\pi^{(l)}=(\pi^{(l)}_i)_{i\in\mathbb{Z}}$ be a $\sigma$-finite
invariant measure corresponding to the family $(P^{(l)}_i)$. Let
$\wtilde
{P}^{(l)}_i$ be the lateral
extension of $P^{(l)}_i$ to $\mathbb{Z}^{\mathbb{Z}}$, that is, under
$\wtilde{P}^{(l)}_i$, $x(0)=i, (x(0),x(1),\ldots)$ is a Markov chain
with transition probabilities $(p^{(l)}_{jk})$ and
$(x(0),x(-1),\ldots)$ is a Markov chain with transition
probabilities $({\pi^{(l)}_k p^{(l)}_{kj}} / \pi^{(l)}_j)$. Assume at
least one (say, the first one) of the Markov chains is null-recurrent
and define
a $\sigma$-finite measure $\mu$ on $S=(\mathbb{Z}^{\mathbb
{Z}})^d$ by
\[
\mu(A_1 \times A_2 \times\cdots\times A_d) = \prod_{l=1}^d
\Biggl(\sum_{i=-\infty}^{\infty}
\pi^{(l)}_i \wtilde{P}^{(l)}_i(A_l)\Biggr),
\]
and observe that $\mu$ is invariant under the $\mathbb{Z}^d$-action
$\{\phi_{(i_1,i_2,\ldots,i_d)}\}_{(i_1,\ldots,i_d)\in\mathbb Z^d}$ on
$S$ defined as the coordinatewise left shift, that is,
%
\begin{equation}\label{defactionex2}
\phi_{(i_1,\ldots,i_d)}\bigl(a^{(1)}, \ldots,
a^{(d)}\bigr)(u_1,\ldots,u_d)=\bigl(a^{(1)}(u_1+i_1),\ldots,
a^{(d)}(u_d+i_d) \bigr)\hspace*{-32pt}
\end{equation}
for all $(a^{(1)}, \ldots,
a^{(d)}) \in S$ and $u_1,\ldots,u_d \in\mathbb{Z}$.

Let $\mathbf{X}=\{X_{(i_1,i_2,\ldots,i_d)}\}_{(i_1,\ldots
,i_d)\in\mathbb Z^d}$
be a stationary S$\alpha$S random field defined by the integral
representation (\ref{repintegralRep}) with $M_\alpha$ being a
S$\alpha
$S random
measure on $S$ with control measure $\mu$ and
\[
f_{(i_1,i_2,\ldots,i_d)}=f \circ
\phi_{(i_1,i_2,\ldots,i_d)},\qquad i_1,i_2,\ldots,i_d \in
\mathbb{Z},
\]
with
\begin{eqnarray}
f\bigl(x^{(1)},x^{(2)},\ldots,x^{(d)}\bigr)=\ind_{\{
x^{(1)}(0)=x^{(2)}(0)=\cdots=
x^{(d)}(0)=0\}},\nonumber\\
&&\eqntext{x^{(1)},x^{(2)},\ldots,x^{(d)} \in
\mathbb{Z}^{\mathbb{Z}}.}
\end{eqnarray}
Clearly, the restriction of (\ref{defactionex2}) to the first
coordinate is a null flow because the first Markov chain is
null-recurrent (see Example 4.1 in~\cite{samorodnitsky05null}) and,
hence, (\ref{defactionex2}) is a null $\mathbb{Z}^d$-action. This
shows, in particular, that $\mathbf{X}$ is weakly mixing. However, if
$d>1$ and some of the Markov chains are positive-recurrent, then the
restriction of $\mu$ in the corresponding coordinate directions are
finite and, hence, by Theorem~\ref{thmsumErgodicity}, $\mathbf{X}$ is
not ergodic along those directions. In this case, the random field
cannot be mixing because it is not mixing in every coordinate
direction. This gives examples of stationary $d$-dimensional ($d>1$)
S$\alpha$S random fields, which are weakly mixing but not mixing. See
Example~4.2 in~\cite{gross93ergodic} for such an example in the $d=1$ case.
\end{Example}
%
\begin{Rem} Correspondingly, we can define $\alpha$-Fr\'echet random
fields and apply Theorem~\ref{thmmaxErgodicity1}. In particular, when
$d>1$, we can obtain an example of an $\alpha$-Fr\'echet random field,
which is weakly mixing but not mixing.
\end{Rem}

%
\begin{appendix}
\section*{Appendix: Proofs of auxiliary results}\label{app}
\subsection{\texorpdfstring{Proof of Lemma \protect\ref{lemmaximal}}{Proof of Lemma 2.2}}
Set
%
\setcounter{equation}{0}
\begin{equation}\label{equ}
u(I(\calG)) \defe\sup_{\nu\in\Lambda(\calG)}\mu(S_\nu).
\end{equation}
Without loss of generality, we assume $\mu(S)<\infty$ (recall that
$\mu
$ is $\sigma$-finite), whence $u(I(\calG))<\infty$. Then, there exists
a sequence of measures $\{\nu_n\}_{n\in\mathbb N}\subset\Lambda
(\calG)$, such that
$u_n\defe\mu(S_{\nu_n})\to u(I(\calG))$ as $n\to\infty$. Set
\[
P_\calG\defe\bcupn S_{\nu_n}.
\]
Clearly, $P_\calG$ is measurable. We show that there exists $\nu
_\calG
\in\Lambda(\calG)$ such that $S_{\nu_\calG} = P_\calG$ and $\mu
(P_\calG
) = u(I(\calG))$. Indeed, we can define on $(S,\calB)$ the measure
%
\begin{equation}\label{eqnucalG}
\nu_\calG(A)\defe\sumn\frac1{2^nu_n}\nu_n(A) \qquad\mbox{for all }
A\in\calB.
\end{equation}
Clearly, $\nu_\calG\in\Lambda(\calG)$, $S_{\nu_\calG} = P_\calG
\mod\mu
$, and $\mu(P_\calG)\leq u(I(\calG))$ by (\ref{equ}). It is also clear
that for all $n\in\mathbb N$, $\nu_n\ll\nu_\calG$ and, hence,
$P_\calG
\supset S_{\nu_n}\mod\mu$. This implies $\mu(P_\calG)\geq u_n$ for all
$n\in\mathbb N$. We have thus shown that $\mu(P_\calG) = u(I(\calG))$.

To complete the proof, we show $P_\calG$ is unique modulo $\mu$-null
sets. Suppose there exist $P_\calG^{(1)}$ and $P_\calG^{(2)}$ such
that $\mu(P_\calG^{(1)}) = \mu(P_\calG^{(2)}) = u(I(\calG))$ and
$\mu
(P_\calG^{(1)}\triangle P_\calG^{(2)}) > 0$. Suppose $\nu^{(1)},\nu
^{(2)}\in\Lambda(\calG)$ are defined as in (\ref{eqnucalG}), so that
$S_{\nu^{(i)}} = P_\calG^{(i)}$ for $i = 1,2$. Clearly, $\nu
^{(1)}+\nu
^{(2)}\in\Lambda(\calG)$. Then, we have $P_\calG^{(1)}\cup
P_\calG^{(2)}\subset I(\calG)$ and $\mu(P_\calG^{(1)}\cup P_\calG
^{(2)})>u(I(\calG
))$, which contradicts (\ref{equ}). The proof is thus
complete.\vspace*{-2pt}

\subsection{\texorpdfstring{Proof of Theorem \protect\ref{thmdecompPN}}{Proof of Theorem 2.3}}
First we introduce some notation. For all transformation $\phi$ on
$\sbm
$, write
\[
\Lambda(\phi)\defe\{\nu\ll\mu\dvtx\nu\mbox{ finite positive measure
on }
S, \nu\circ\phi\inv= \nu\}.
\]
We need the following lemma.\vspace*{-3pt}
%
\begin{Lem}\label{lem2}
Suppose $\phi$ is an arbitrary invertible, bimeasurable and nonsingular
transformation on $\sbm$. Then $\mu(\phi\inv(S_\nu)\triangle S_\nu
) = 0$, for all $\nu\in\Lambda(\phi)$.\vspace*{-3pt}
\end{Lem}
\begin{pf}
First, we show for all $\nu\in\Lambda(\phi)$, $\mu(\phi\inv
(S_\nu
)\triangle S_\nu) = 0$. If not, then set $E_0 \defe\phi\inv(S_\nu
)\setminus S_\nu$, $F_0 = \phi(E_0)$ and suppose $\mu(E_0)>0$. Since
$\phi$ is nonsingular, $\mu(F_0) >0$. Note that $F_0\subset S_\nu$ and
$\mu\sim\nu$ on $S_\nu$, whence $\nu(F_0)>0$. Note also that $\nu
(S_\nu
^c) = 0$ and $\nu\circ\phi\inv= \nu$ imply $\nu(F_0) = \nu\circ
\phi\inv
(F_0) = \nu(E_0) \leq\nu(S_\nu^c) = 0$. This contradicts $\nu(F_0)>0$.
We have thus shown that $\mu(\phi\inv(S_\nu)\setminus S_\nu) = 0$.

Next, we show that $\mu(S_\nu\setminus\phi\inv(S_\nu)) = 0$. Indeed,
setting $E_1 \defe S_\nu\setminus\phi\inv(S_\nu)$, we have $\nu
(S_\nu)
= \nu(E_1) + \nu(\phi\inv(S_\nu)\cap S_\nu)$. At the same time,
$\nu
(S_\nu) = \nu\circ\phi\inv(S_\nu) = \nu(\phi\inv(S_\nu)\cap
S_\nu) + \nu
(E_0)$, where $E_0\defe\phi\inv(S_\nu)\setminus S_\nu$. Since $\nu
(E_0) = 0$ as shown in the first part of the proof, the two equations
above imply $\nu(E_1) = 0$, since $\nu$ is finite. Finally, by the fact
that $\nu\sim\mu$ on $S_\nu$, we have $\mu(S_\nu\setminus\phi
\inv(S_\nu
)) \equiv\mu(E_1) = 0$.\vspace*{-2pt}
\end{pf}

Now we prove Theorem~\ref{thmdecompPN}.\vspace*{-2pt}

\begin{longlist}[(iii)]
\item[(i)]
Fix $\phi\in\calG$. Note that by Lemma~\ref{lemmaximal},
there exists $\nu_{\calG}\in\Lambda(\phi)\subset I(\calG)$ such that
$S_{\nu_{\calG}} = P_\calG$. Then, by Lemma~\ref{lem2}, $\mu(\phi
\inv
(P_\calG)\triangle P_\calG) = 0$. By the fact that all $\phi\in
\calG$
are invertible, we have that $\phi\inv(N_\calG)^c = \phi\inv
(N_\calG
^c)$ and by the identity $A\triangle B = A^c\triangle B^c$, we have
$\mu
(\phi\inv(N_\calG)\triangle N_\calG) = 0$. The previous argument is
valid for all $\phi\in\calG$.\vadjust{\goodbreak}

\item[(ii)] Consider $L^1(P_\calG,\calB\cap P_\calG,\mu|_{P_\calG
})$, where $\calB\cap P_\calG\defe\{A\cap P_\calG\dvtx A\in\calB\}$ and
$\mu|_{P_\calG}$ is the restriction of $\mu$ tn $\calB\cap P_\calG
$. Define
%
\begin{eqnarray}\label{eqdualV1}
\widetilde\phi f(s)\equiv[\widetilde\phi(f)](s)\defe
{\frac{\d(\mu\circ\phi\inv)}{\d\mu}}(s)f\circ\phi\inv(s)\ind
_{P_\calG
\cap\phi
(P_\calG)}(s)\nonumber\\[-8pt]\\[-8pt]
&&\eqntext{\mbox{for all } f\in L^1(P_\calG,\mu|_{P_\calG}).}
\end{eqnarray}
In this way, the mapping $\widetilde\phi$ is a restricted version of
$\widehat\phi$ on $L^1(P_\calG,\mu|_{P_\calG})$ in the sense that
%
\begin{equation}\label{eqwidetildephi}
\widetilde\phi f = \widehat\phi f,\qquad\mu|_{P_\calG}\mbox{-a.e.}\mbox{ for
all }
f\in L^1(P_\calG,\mu|_{P_\calG})\subset L^1(S,\mu).
\end{equation}

Recall\vspace*{1pt} that by Lemma~\ref{lemmaximal} there exists $\nu\in\Lambda
(\calG
)$ such that $\widehat\phi(\d\nu/\d\mu) = \d\nu/\d\mu$ for
all $\phi\in
\calG$ and $\supp(\nu) = P_\calG$. Whence, for $\widetilde\nu
\defe\nu
|_{P_\calG}$, we have $\widetilde\phi(\d\widetilde\nu/\d\mu
|_{P_\calG})
= \d\widetilde\nu/\d\mu|_{P_\calG}$ for all $\phi\in\calG$ and
$\widetilde\nu\sim\mu|_{P_\calG}$.
Note that all locally compact Abelian groups are amenable (see, e.g.,
Example 1.1.5(c) in~\cite{runde02lectures}). Thus, Theorem 1 [parts (1)
and~(8)] in~\cite{takahashi71invariant} applied to $\widetilde G$ and
$f$ implies that
\[
\sumn\widetilde\phi_{u_n}f(s) = \infty,
\qquad\mu|_{P_\calG}\mbox{-a.e.}\mbox
{ for
all } \{\widetilde\phi_{u_n}\}_{n\in\mathbb N}\subset\widetilde
\calG,
\]
which, by (\ref{eqwidetildephi}), is equivalent to (\ref
{eqPG}).

\item[(iii)] Similarly, as in (ii), restrict $\calG$ to $L^1(N_\calG
,\calB\cap N_\calG,\mu|_{N_\calG})$ and apply Theorem~2 [parts (1) and
(8)] in~\cite{takahashi71invariant}.
\end{longlist}

\subsection{\texorpdfstring{Proof of Theorem \protect\ref{thmweaklyWandering}}{Proof of Theorem 2.4}}
We only sketch the proof of this result.

\begin{longlist}[(ii)]
\item[(i)]
We apply\vspace*{1pt} Theorem 1 [parts (1) and (6)] in \cite
{takahashi71invariant}.
Recall that the adjoint operator of $\widehat\phi$, $\widehat\phi
^*\dvtx(L^1)^*\to(L^1)^*$ [$(L^1)^* = L^\infty$]
is such that for all $f\in\lone\smu$ and $h\in L^\infty\smu$,
\[
\int_Sf(s)[\widehat\phi^*(h)](s)\mu(\d s) = \int_S[\widehat\phi
(f)](s)h(s)\mu(\d s).
\]
The last integral equals
\[
\int_S\frac{\d(\mu\circ\phi\inv)}{\d\mu}(s)f\circ\phi\inv
(s)h\circ
\phi\circ\phi\inv
(s)\mu(\d s) = \int_Sf(s)h\circ\phi(s)\mu(\d s),
\]
whence $[\widehat\phi^*(h)](s) = h\circ\phi(s),\mu$-a.e. Thus,
if $W$
is a weakly wandering set w.r.t.~$\calG$, we have
\[
\sumn\widehat\phi_{t_n}^*\ind_W(s)<2\qquad\mbox{for some } \{\phi
_{t_n}\}
_{n\in\mathbb N}\subset\calG.
\]
Now, part (6) of Theorem 1 in~\cite{takahashi71invariant} is equivalent
to the nonexistence of a weakly wandering set of positive measure.

\item[(ii)] The proof is similar to the proof of Proposition 1.4.7\vadjust{\goodbreak}
in~\cite{aaronson97introduction}.
\end{longlist}


\subsection{\texorpdfstring{Proof of Proposition \protect\ref{propgross}}{Proof of Proposition 4.2}}
We first need the following lemma.
%
\begin{Lem}\label{lemgross}
Assume $\{X_t\}_{t\in\mathbb T^d}$ is a stationary S$\alpha$S random
field with spectral
representation $\{f_t\}_{t\in\mathbb T^d}\subset\la\sbm$, $\alpha
\in(0,2)$. Then,
$\{X_t\}_{t\in\mathbb T^d}$ is weakly mixing, if and only if, there
exists a density one
set $D\subset\mathbb T^d$, such that
%
\begin{eqnarray}\label{eqgrossPQ}
\limn\mu\Biggl\{s\dvtx\Biggl|\sum_{j=1}^p\beta_jf_{\tau_j}(s)\Biggr|\in
K, \Biggl|\sum_{k=1}^q\gamma_kf_{t_k+t^*_n}(s)\Biggr|> \varepsilon
\Biggr\} = 0\nonumber\\[3pt]\\[-20pt]
&&\eqntext{\mbox{for all } p,q\in\mathbb N,\beta_j,\gamma_k\in\mathbb R,\tau_j,
t_k\in\mathbb T^d,}\\
\eqntext{\mbox{compact } K\subset\mathbb R\setminus\{0\}, \varepsilon>0
\mbox{ and }
\{t^*_n\}_{n\in\mathbb N}\in\calT_D.}
\end{eqnarray}
\end{Lem}
\begin{pf} It transpires from the proofs in~\cite{maruyama70infinitely}
that a stationary process $\{X_t\}_{t\in\mathbb T^d}$ is weakly mixing
if and only if
there exists a density one set $D\subset\mathbb T^d$ such that
%
\begin{eqnarray}\label{eqmaruyama70}
&&\limn\esp{\Biggl[\exp\Biggl(i\sumjp\beta_jX_{\tau
_j}\Biggr)\exp\Biggl(i\sum
_{k=1}^q\gamma_kX_{t_k+t^*_n}\Biggr)\Biggr]} \nonumber\\
&&\qquad = \esp\exp\Biggl(i\sumjp\beta_jX_{\tau_j}\Biggr)\esp\exp\Biggl(i\sum
_{k=1}^q\gamma_kX_{t_k}\Biggr)\\
\eqntext{\mbox{for all } p,q\in\mathbb N, \beta_j,\gamma_k\in\mathbb R,
\tau
_j,t_k\in\mathbb T\mbox{ and }\{t^*_n\}_{n\in\mathbb
N}\in\calT_D.}
\end{eqnarray}
See the following remark on the equivalence of (\ref{eqgrossPQ}) and
(\ref{eqmaruyama70}).
\end{pf}
%
\begin{Rem} \label{remarkmistakegross}
In the one-dimensional case, to show that (\ref{eqgross94}) is
equivalent to the weak mixing of the process, Gross~\cite{gross94some}
proved that (\ref{eqgross94}) is equivalent to the following weaker
condition (\ref{eqmaruyama70}) (Theorem 2.7 in~\cite{gross94some}):
%
\begin{eqnarray}\label{eqmaruyama94}
\limn\esp[\exp(i\theta_1X_{0})\exp(i\theta_2 X_{t_n})] =
\esp\exp
(i\theta_1 X_0)\esp\exp(i\theta_2 X_0)\nonumber\\[-8pt]\\[-8pt]
\eqntext{\mbox{for all }\theta_1,\theta_2\in\mathbb R,\{t_n\}_{n\in\mathbb
N}\in\calT_D.}
\end{eqnarray}
The equivalence of (\ref{eqmaruyama70}) and (\ref{eqmaruyama94}),
however, seems nontrivial and yet not mentioned in~\cite{gross94some}.
Nevertheless, parallel to the proof of Theorem 2.7 in \cite
{gross94some}, we can prove Lemma~\ref{lemgross}.
\end{Rem}

To show Proposition~\ref{propgross}, it suffices to prove the
following lemma.
%
\begin{Lem}\label{lemgap}
Assume $\alpha\in(0,2)$ and $\{X_t\}_{t\in\mathbb T^d}$ is a
stationary S$\alpha$S
process with spectral representation $\{f_t\}_{t\in\mathbb T^d}\subset
\la\sbm$.
Then (\ref{eqgrossPQ}) is true if and only if (\ref{eqgross94}) is
true.
\end{Lem}
\begin{pf}
Clearly, (\ref{eqgrossPQ}) implies (\ref{eqgross94}). Now suppose
that (\ref{eqgross94}) is true. We will show (\ref{eqgrossPQ}). For
any $p,q\in\mathbb N$ and $\tau_j,t_k\in\mathbb T^d$, write
%
\begin{equation}\label{eqgphq}
g_p(s)\defe\sumjp\beta_jf_{\tau_j}(s)\quad\mbox{and}\quad h_q(s) \defe\sum
_{k=1}^q\gamma_kf_{t_k}(s).
\end{equation}
We will prove (\ref{eqgrossPQ}) by induction on $(p,q)$. By (\ref
{eqgross94}), we have that (\ref{eqgrossPQ}) holds for $(p,q) = (1,1)$.

\begin{longlist}[(ii)]
\item[(i)]
Suppose for fixed $(p,q)$ (\ref{eqgrossPQ}) holds, then
we will show that (\ref{eqgrossPQ}) holds for $(p+1,q)$. If not, then
there exists $\{t^*_n\}_{n\in\mathbb N}\in\calT_D$ such that for
some compact
$K\subset
\mathbb R\setminus\{0\}$ and $\delta>0$, we have $\mu(E_n)\geq
\delta$ with
\[
E_n\defe\{s\dvtx|g_p(s)+\beta_{p+1}f_{\tau_{p+1}}(s)|\in K
,|U_{t^*_n}h_q(s)|>\varepsilon\}.
\]
Here for all $t\in\mathbb T^d$, $U_t(\sum_{k=1}^q\gamma
_kf_{t_k})(s) \defe\sum_{k=1}^q\gamma_kf_{t_k+t}(s)$.

Without loss of generality, we can assume $K\subset(0,\infty)$. Then,
since $K$ is compact, there exists $0<d_K<M$ such that $K\subset[d_K,M]$.
Since $f_{\tau_1},\ldots,\break f_{\tau_{p+1}}\in\la(S,\mu)$, we can also
choose $M$ to be large enough so that $\mu(E_M^0)\leq\delta/2$, where
\[
E_M^0\defe\{s\dvtx|g_p(s)|>M \mbox{ or } |\beta_{p+1}f_{\tau
_{p+1}}(s)|>M\}.
\]
Then, we claim that for each $n$, either of the two sets
\[
E_{n}^p\defe\biggl\{s\dvtx|g_p(s)|\in\biggl[\frac{d_K}2,M
\biggr],|U_{t^*_n}h_q(s)|>\varepsilon\biggr\}
\]
and
\[
E_n^{p+1}\defe\biggl\{s\dvtx|\beta_{p+1}f_{\tau_{p+1}}(s)|\in
\biggl[\frac{d_K}2,M\biggr],|U_{t^*_n}h_g(s)|>\varepsilon\biggr\}
\]
has measure larger than $\delta/4$. Otherwise, observe that
$E_n\subset
E_{n}^p\cup E_{n}^{p+1}\cup E_M^0$,
which implies that $\mu(E_n)<\delta$, a contradiction.

It then follows that either $\{E^p_n\}_{n\in\mathbb N}$ or $\{
E^{p+1}_n\}_{n\in\mathbb N}$ will have
a subsequence with measures larger than $\delta/4$. Namely, there
exists $\{t^*_{n_k}\}_{k\in\mathbb N}\in\calT_D$ such that
\[
\mu(E_{n_k}^p)\geq\frac\delta4\qquad\mbox{for all } k\in\mathbb N
\quad\mbox{or}\quad
\mu
(E_{n_k}^{p+1})\geq\frac\delta4\qquad\mbox{for all } k\in\mathbb N.
\]
But the first case contradicts the assumption that (\ref{eqgrossPQ})
holds for $(p,q)$ and the second case contradicts (\ref{eqgross94}).
We have thus shown that (\ref{eqgrossPQ}) holds for $(p+1,q)$.

\item[(ii)]
Next, suppose (\ref{eqgrossPQ}) holds for $(p,q)$ and
we show that it holds for $(p$, \mbox{$q+1)$}. If not, then there exists a
compact $K\subset\mathbb R\setminus\{0\}$ such that
\[
\mu\{s\dvtx|g_p(s)|\in K,|U_{t^*_n}(h_q+\gamma
_{q+1}f_{t_{q+1}})(s)|>\varepsilon\}\nrightarrow0
\qquad\mbox{as } n\to\infty
.
\]
Then, by a similar argument as in part (i), one can show that for all
$\varepsilon>0$, there exists $\{t^*_n\}_{n\in\mathbb N}\in\calT_D$
and $\delta>0$ such
that we have either
\[
\mu\biggl\{s\dvtx|g_p(s)|\in K, |U_{t^*_{n}}h_q(s)|>\frac\varepsilon2\biggr\}
\geq
\delta>0
\]
or
\[
\mu\biggl\{s\dvtx|g_p(s)|\in K, |\gamma_{q+1}f_{t_{q+1}+t^*_n}(s)|>\frac
\varepsilon2\biggr\}\geq\delta>0.
\]
Both cases lead to contradictions. We have thus shown that (\ref
{eqgrossPQ}) holds for $(p,q+1)$. The proof is thus complete.\qed
\end{longlist}
\noqed\end{pf}
\end{appendix}

\section*{Acknowledgments}

The authors are thankful to Jan
Rosi\'nski for suggesting the problem of equivalence of ergodicity and
weak mixing for the max-stable case, and to Yimin Xiao for a number of
useful discussions on the properties of the local times of a Brownian
sheet. The authors are also thankful to the anonymous referees for
their detailed comments.



\printaddresses


\begin{thebibliography}{43}

\bibitem{aaronson97introduction}
\begin{bbook}[mr]
\bauthor{\bsnm{Aaronson},~\bfnm{Jon}\binits{J.}}
(\byear{1997}).
\btitle{An Introduction to Infinite Ergodic Theory}.
\bseries{Mathematical Surveys and Monographs}
\bvolume{50}.
\bpublisher{Amer. Math. Soc.}, \baddress{Providence, RI}.
\bid{mr={1450400}}
\bptok{imsref}%
\end{bbook}
\endbibitem

\bibitem{cambanis87ergodic}
\begin{barticle}[mr]
\bauthor{\bsnm{Cambanis},~\bfnm{Stamatis}\binits{S.}},
  \bauthor{\bsnm{Hardin},~\bfnm{Clyde~D.}\binits{C.~D.} \bsuffix{Jr.}} \AND
  \bauthor{\bsnm{Weron},~\bfnm{Aleksander}\binits{A.}}
(\byear{1987}).
\btitle{Ergodic properties of stationary stable processes}.
\bjournal{Stochastic Process. Appl.}
\bvolume{24}
\bpages{1--18}.
\bid{doi={10.1016/0304-4149(87)90024-X}, issn={0304-4149}, mr={0883599}}
\bptok{imsref}%
\end{barticle}
\endbibitem

\bibitem{cohen06random}
\begin{barticle}[mr]
\bauthor{\bsnm{Cohen},~\bfnm{Serge}\binits{S.}} \AND
  \bauthor{\bsnm{Samorodnitsky},~\bfnm{Gennady}\binits{G.}}
(\byear{2006}).
\btitle{Random rewards, fractional {B}rownian local times and stable
  self-similar processes}.
\bjournal{Ann. Appl. Probab.}
\bvolume{16}
\bpages{1432--1461}.
\bid{doi={10.1214/105051606000000277}, issn={1050-5164}, mr={2260069}}
\bptok{imsref}%
\end{barticle}
\endbibitem

\bibitem{dehaan84spectral}
\begin{barticle}[mr]
\bauthor{\bparticle{de} \bsnm{Haan},~\bfnm{L.}\binits{L.}}
(\byear{1984}).
\btitle{A spectral representation for max-stable processes}.
\bjournal{Ann. Probab.}
\bvolume{12}
\bpages{1194--1204}.
\bid{issn={0091-1798}, mr={0757776}}
\bptok{imsref}%
\end{barticle}
\endbibitem

\bibitem{dehaan86stationary}
\begin{barticle}[mr]
\bauthor{\bparticle{de} \bsnm{Haan},~\bfnm{L.}\binits{L.}} \AND
  \bauthor{\bsnm{Pickands},~\bfnm{J.}\binits{J.} \bsuffix{III}}
(\byear{1986}).
\btitle{Stationary min-stable stochastic processes}.
\bjournal{Probab. Theory Related Fields}
\bvolume{72}
\bpages{477--492}.
\bid{doi={10.1007/BF00344716}, issn={0178-8051}, mr={0847381}}
\bptok{imsref}%
\end{barticle}
\endbibitem

\bibitem{ehm81sample}
\begin{barticle}[mr]
\bauthor{\bsnm{Ehm},~\bfnm{W.}\binits{W.}}
(\byear{1981}).
\btitle{Sample function properties of multiparameter stable processes}.
\bjournal{Z.~Wahrsch. Verw. Gebiete}
\bvolume{56}
\bpages{195--228}.
\bid{doi={10.1007/BF00535741}, issn={0044-3719}, mr={0618272}}
\bptok{imsref}%
\end{barticle}
\endbibitem

\bibitem{gross94some}
\begin{barticle}[mr]
\bauthor{\bsnm{Gross},~\bfnm{Aaron}\binits{A.}}
(\byear{1994}).
\btitle{Some mixing conditions for stationary symmetric stable stochastic
  processes}.
\bjournal{Stochastic Process. Appl.}
\bvolume{51}
\bpages{277--295}.
\bid{doi={10.1016/0304-4149(94)90046-9}, issn={0304-4149}, mr={1288293}}
\bptok{imsref}%
\end{barticle}
\endbibitem

\bibitem{gross93ergodic}
\begin{barticle}[mr]
\bauthor{\bsnm{Gross},~\bfnm{Aaron}\binits{A.}} \AND
  \bauthor{\bsnm{Robertson},~\bfnm{James~B.}\binits{J.~B.}}
(\byear{1993}).
\btitle{Ergodic properties of random measures on stationary sequences of sets}.
\bjournal{Stochastic Process. Appl.}
\bvolume{46}
\bpages{249--265}.
\bid{doi={10.1016/0304-4149(93)90006-P}, issn={0304-4149}, mr={1226411}}
\bptok{imsref}%
\end{barticle}
\endbibitem

\bibitem{kabluchko09spectral}
\begin{barticle}[mr]
\bauthor{\bsnm{Kabluchko},~\bfnm{Zakhar}\binits{Z.}}
(\byear{2009}).
\btitle{Spectral representations of sum- and max-stable processes}.
\bjournal{Extremes}
\bvolume{12}
\bpages{401--424}.
\bid{doi={10.1007/s10687-009-0083-9}, issn={1386-1999}, mr={2562988}}
\bptok{imsref}%
\end{barticle}
\endbibitem

\bibitem{kabluchko10ergodic}
\begin{barticle}[mr]
\bauthor{\bsnm{Kabluchko},~\bfnm{Zakhar}\binits{Z.}} \AND
  \bauthor{\bsnm{Schlather},~\bfnm{Martin}\binits{M.}}
(\byear{2010}).
\btitle{Ergodic properties of max-infinitely divisible processes}.
\bjournal{Stochastic Process. Appl.}
\bvolume{120}
\bpages{281--295}.
\bid{doi={10.1016/j.spa.2009.12.002}, issn={0304-4149}, mr={2584894}}
\bptok{imsref}%
\end{barticle}
\endbibitem

\bibitem{kolodynski03group}
\begin{barticle}[mr]
\bauthor{\bsnm{Kolody{\'n}ski},~\bfnm{S{\l}awomir}\binits{S.}} \AND
  \bauthor{\bsnm{Rosi{\'n}ski},~\bfnm{Jan}\binits{J.}}
(\byear{2003}).
\btitle{Group self-similar stable processes in {${\Bbb R}\sp d$}}.
\bjournal{J. Theoret. Probab.}
\bvolume{16}
\bpages{855--876}.
\bid{doi={10.1023/B:JOTP.0000011997.14357.fd}, issn={0894-9840}, mr={2033189}}
\bptok{imsref}%
\end{barticle}
\endbibitem

\bibitem{krengel67classification}
\begin{bincollection}[mr]
\bauthor{\bsnm{Krengel},~\bfnm{Ulrich}\binits{U.}}
(\byear{1967}).
\btitle{Classification of states for operators}.
In \bbooktitle{Proc. {F}ifth {B}erkeley {S}ympos. {M}ath. {S}tatist. and
  {P}robability ({B}erkeley, {C}alif., 1965/66), {V}ol. {II}: {C}ontributions
  to {P}robability {T}heory, {P}art 2}
\bpages{415--429}.
\bpublisher{Univ. California Press}, \baddress{Berkeley, CA}.
\bid{mr={0241601}}
\bptok{imsref}%
\end{bincollection}
\endbibitem

\bibitem{krengel85ergodic}
\begin{bbook}[mr]
\bauthor{\bsnm{Krengel},~\bfnm{Ulrich}\binits{U.}}
(\byear{1985}).
\btitle{Ergodic Theorems}.
\bseries{de Gruyter Studies in Mathematics}
\bvolume{6}.
\bpublisher{de Gruyter}, \baddress{Berlin}.
\bnote{With a supplement by Antoine Brunel}.
\bid{doi={10.1515/9783110844641}, mr={0797411}}
\bptok{imsref}%
\end{bbook}
\endbibitem

\bibitem{maruyama70infinitely}
\begin{barticle}[mr]
\bauthor{\bsnm{Maruyama},~\bfnm{G.}\binits{G.}}
(\byear{1970}).
\btitle{Infinitely divisible processes}.
\bjournal{Teor. Verojatnost. i Primenen.}
\bvolume{15}
\bpages{3--23}.
\bid{issn={0040-361X}, mr={0285046}}
\bptok{imsref}%
\end{barticle}
\endbibitem

\bibitem{neveu67existence}
\begin{bincollection}[mr]
\bauthor{\bsnm{Neveu},~\bfnm{Jacques}\binits{J.}}
(\byear{1967}).
\btitle{Existence of bounded invariant measures in ergodic theory}.
In \bbooktitle{Proc. {F}ifth {B}erkeley {S}ympos. {M}ath. {S}tatist. and
  {P}robability ({B}erkeley, {C}alif., 1965/66), {V}ol. {II}: {C}ontributions
  to {P}robability {T}heory, {P}art 2}
\bpages{461--472}.
\bpublisher{Univ. California Press}, \baddress{Berkeley, CA}.
\bid{mr={0212161}}
\bptok{imsref}%
\end{bincollection}
\endbibitem

\bibitem{petersen83ergodic}
\begin{bbook}[mr]
\bauthor{\bsnm{Petersen},~\bfnm{Karl}\binits{K.}}
(\byear{1983}).
\btitle{Ergodic Theory}.
\bseries{Cambridge Studies in Advanced Mathematics}
\bvolume{2}.
\bpublisher{Cambridge Univ. Press}, \baddress{Cambridge}.
\bid{mr={0833286}}
\bptok{imsref}%
\end{bbook}
\endbibitem

\bibitem{pipiras07nonminimal}
\begin{barticle}[mr]
\bauthor{\bsnm{Pipiras},~\bfnm{Vladas}\binits{V.}}
(\byear{2007}).
\btitle{Nonminimal sets, their projections and integral representations of
  stable processes}.
\bjournal{Stochastic Process. Appl.}
\bvolume{117}
\bpages{1285--1302}.
\bid{doi={10.1016/j.spa.2007.01.002}, issn={0304-4149}, mr={2343940}}
\bptok{imsref}%
\end{barticle}
\endbibitem

\bibitem{pipiras04stable}
\begin{barticle}[mr]
\bauthor{\bsnm{Pipiras},~\bfnm{Vladas}\binits{V.}} \AND
  \bauthor{\bsnm{Taqqu},~\bfnm{Murad~S.}\binits{M.~S.}}
(\byear{2004}).
\btitle{Stable stationary processes related to cyclic flows}.
\bjournal{Ann. Probab.}
\bvolume{32}
\bpages{2222--2260}.
\bid{doi={10.1214/009117904000000108}, issn={0091-1798}, mr={2073190}}
\bptok{imsref}%
\end{barticle}
\endbibitem

\bibitem{podgorski92note}
\begin{barticle}[mr]
\bauthor{\bsnm{Podg{\'o}rski},~\bfnm{Krzysztof}\binits{K.}}
(\byear{1992}).
\btitle{A note on ergodic symmetric stable processes}.
\bjournal{Stochastic Process. Appl.}
\bvolume{43}
\bpages{355--362}.
\bid{doi={10.1016/0304-4149(92)90068-2}, issn={0304-4149}, mr={1191157}}
\bptok{imsref}%
\end{barticle}
\endbibitem

\bibitem{rosinski95structure}
\begin{barticle}[mr]
\bauthor{\bsnm{Rosi{\'n}ski},~\bfnm{Jan}\binits{J.}}
(\byear{1995}).
\btitle{On the structure of stationary stable processes}.
\bjournal{Ann. Probab.}
\bvolume{23}
\bpages{1163--1187}.
\bid{issn={0091-1798}, mr={1349166}}
\bptok{imsref}%
\end{barticle}
\endbibitem

\bibitem{rosinski00decomposition}
\begin{barticle}[mr]
\bauthor{\bsnm{Rosi{\'n}ski},~\bfnm{Jan}\binits{J.}}
(\byear{2000}).
\btitle{Decomposition of stationary {$\alpha$}-stable random fields}.
\bjournal{Ann. Probab.}
\bvolume{28}
\bpages{1797--1813}.
\bid{doi={10.1214/aop/1019160508}, issn={0091-1798}, mr={1813849}}
\bptok{imsref}%
\end{barticle}
\endbibitem

\bibitem{rosinski06minimal}
\begin{barticle}[mr]
\bauthor{\bsnm{Rosi{\'n}ski},~\bfnm{Jan}\binits{J.}}
(\byear{2006}).
\btitle{Minimal integral representations of stable processes}.
\bjournal{Probab. Math. Statist.}
\bvolume{26}
\bpages{121--142}.
\bid{issn={0208-4147}, mr={2301892}}
\bptok{imsref}%
\end{barticle}
\endbibitem

\bibitem{rosinski96classes}
\begin{barticle}[mr]
\bauthor{\bsnm{Rosi{\'n}ski},~\bfnm{Jan}\binits{J.}} \AND
  \bauthor{\bsnm{Samorodnitsky},~\bfnm{Gennady}\binits{G.}}
(\byear{1996}).
\btitle{Classes of mixing stable processes}.
\bjournal{Bernoulli}
\bvolume{2}
\bpages{365--377}.
\bid{doi={10.2307/3318419}, issn={1350-7265}, mr={1440274}}
\bptok{imsref}%
\end{barticle}
\endbibitem

\bibitem{rosinski96simple}
\begin{barticle}[mr]
\bauthor{\bsnm{Rosi{\'n}ski},~\bfnm{Jan}\binits{J.}} \AND \bauthor{\bsnm{{\.
  Z}ak},~\bfnm{Tomasz}\binits{T.}}
(\byear{1996}).
\btitle{Simple conditions for mixing of infinitely divisible processes}.
\bjournal{Stochastic Process. Appl.}
\bvolume{61}
\bpages{277--288}.
\bid{doi={10.1016/0304-4149(95)00083-6}, issn={0304-4149}, mr={1386177}}
\bptok{imsref}%
\end{barticle}
\endbibitem

\bibitem{rosinski97equivalence}
\begin{barticle}[mr]
\bauthor{\bsnm{Rosi{\'n}ski},~\bfnm{Jan}\binits{J.}} \AND \bauthor{\bsnm{{\.
  Z}ak},~\bfnm{Tomasz}\binits{T.}}
(\byear{1997}).
\btitle{The equivalence of ergodicity of weak mixing for infinitely divisible
  processes}.
\bjournal{J. Theoret. Probab.}
\bvolume{10}
\bpages{73--86}.
\bid{doi={10.1023/A:1022690230759}, issn={0894-9840}, mr={1432616}}
\bptok{imsref}%
\end{barticle}
\endbibitem

\bibitem{roy07ergodic}
\begin{barticle}[mr]
\bauthor{\bsnm{Roy},~\bfnm{Emmanuel}\binits{E.}}
(\byear{2007}).
\btitle{Ergodic properties of {P}oissonian {ID} processes}.
\bjournal{Ann. Probab.}
\bvolume{35}
\bpages{551--576}.
\bid{doi={10.1214/009117906000000692}, issn={0091-1798}, mr={2308588}}
\bptok{imsref}%
\end{barticle}
\endbibitem

\bibitem{roy09poisson}
\begin{barticle}[mr]
\bauthor{\bsnm{Roy},~\bfnm{Emmanuel}\binits{E.}}
(\byear{2009}).
\btitle{Poisson suspensions and infinite ergodic theory}.
\bjournal{Ergodic Theory Dynam. Systems}
\bvolume{29}
\bpages{667--683}.
\bid{doi={10.1017/S0143385708080279}, issn={0143-3857}, mr={2486789}}
\bptok{imsref}%
\end{barticle}
\endbibitem

\bibitem{roy10ergodic}
\begin{barticle}[mr]
\bauthor{\bsnm{Roy},~\bfnm{Parthanil}\binits{P.}}
(\byear{2010}).
\btitle{Ergodic theory, Abelian groups and point processes induced by stable
  random fields}.
\bjournal{Ann. Probab.}
\bvolume{38}
\bpages{770--793}.
\bid{doi={10.1214/09-AOP495}, issn={0091-1798}, mr={2642891}}
\bptok{imsref}%
\end{barticle}
\endbibitem

\bibitem{roy10nonsingular}
\begin{barticle}[mr]
\bauthor{\bsnm{Roy},~\bfnm{Parthanil}\binits{P.}}
(\byear{2010}).
\btitle{Nonsingular group actions and stationary {$S\alpha S$} random fields}.
\bjournal{Proc. Amer. Math. Soc.}
\bvolume{138}
\bpages{2195--2202}.
\bid{doi={10.1090/S0002-9939-10-10250-0}, issn={0002-9939}, mr={2596059}}
\bptok{imsref}%
\end{barticle}
\endbibitem

\bibitem{roy08stationary}
\begin{barticle}[mr]
\bauthor{\bsnm{Roy},~\bfnm{Parthanil}\binits{P.}} \AND
  \bauthor{\bsnm{Samorodnitsky},~\bfnm{Gennady}\binits{G.}}
(\byear{2008}).
\btitle{Stationary symmetric {$\alpha$}-stable discrete parameter random
  fields}.
\bjournal{J. Theoret. Probab.}
\bvolume{21}
\bpages{212--233}.
\bid{doi={10.1007/s10959-007-0107-9}, issn={0894-9840}, mr={2384479}}
\bptok{imsref}%
\end{barticle}
\endbibitem

\bibitem{runde02lectures}
\begin{bbook}[mr]
\bauthor{\bsnm{Runde},~\bfnm{Volker}\binits{V.}}
(\byear{2002}).
\btitle{Lectures on Amenability}.
\bseries{Lecture Notes in Math.}
\bvolume{1774}.
\bpublisher{Springer}, \baddress{Berlin}.
\bid{doi={10.1007/b82937}, mr={1874893}}
\bptok{imsref}%
\end{bbook}
\endbibitem

\bibitem{samorodnitsky04extreme}
\begin{barticle}[mr]
\bauthor{\bsnm{Samorodnitsky},~\bfnm{Gennady}\binits{G.}}
(\byear{2004}).
\btitle{Extreme value theory, ergodic theory and the boundary between short
  memory and long memory for stationary stable processes}.
\bjournal{Ann. Probab.}
\bvolume{32}
\bpages{1438--1468}.
\bid{doi={10.1214/009117904000000261}, issn={0091-1798}, mr={2060304}}
\bptok{imsref}%
\end{barticle}
\endbibitem

\bibitem{samorodnitsky04maxima}
\begin{barticle}[mr]
\bauthor{\bsnm{Samorodnitsky},~\bfnm{Gennady}\binits{G.}}
(\byear{2004}).
\btitle{Maxima of continuous-time stationary stable processes}.
\bjournal{Adv. in Appl. Probab.}
\bvolume{36}
\bpages{805--823}.
\bid{doi={10.1239/aap/1093962235}, issn={0001-8678}, mr={2079915}}
\bptok{imsref}%
\end{barticle}
\endbibitem

\bibitem{samorodnitsky05null}
\begin{barticle}[mr]
\bauthor{\bsnm{Samorodnitsky},~\bfnm{Gennady}\binits{G.}}
(\byear{2005}).
\btitle{Null flows, positive flows and the structure of stationary symmetric
  stable processes}.
\bjournal{Ann. Probab.}
\bvolume{33}
\bpages{1782--1803}.
\bid{doi={10.1214/009117905000000305}, issn={0091-1798}, mr={2165579}}
\bptok{imsref}%
\end{barticle}
\endbibitem

\bibitem{samorodnitsky94stable}
\begin{bbook}[mr]
\bauthor{\bsnm{Samorodnitsky},~\bfnm{Gennady}\binits{G.}} \AND
  \bauthor{\bsnm{Taqqu},~\bfnm{Murad~S.}\binits{M.~S.}}
(\byear{1994}).
\btitle{Stable Non-{G}aussian Random Processes:
Stochastic Models with Infinite Variance}.
\bpublisher{Chapman \& Hall}, \baddress{New York}.
\bid{mr={1280932}}
\bptok{imsref}%
\end{bbook}
\endbibitem

\bibitem{stoev08ergodicity}
\begin{barticle}[mr]
\bauthor{\bsnm{Stoev},~\bfnm{Stilian~A.}\binits{S.~A.}}
(\byear{2008}).
\btitle{On the ergodicity and mixing of max-stable processes}.
\bjournal{Stochastic Process. Appl.}
\bvolume{118}
\bpages{1679--1705}.
\bid{doi={10.1016/j.spa.2007.10.013}, issn={0304-4149}, mr={2442375}}
\bptok{imsref}%
\end{barticle}
\endbibitem

\bibitem{stoev06extremal}
\begin{barticle}[mr]
\bauthor{\bsnm{Stoev},~\bfnm{Stilian~A.}\binits{S.~A.}} \AND
  \bauthor{\bsnm{Taqqu},~\bfnm{Murad~S.}\binits{M.~S.}}
(\byear{2005}).
\btitle{Extremal stochastic integrals: A parallel between max-stable processes
  and {$\alpha$}-stable processes}.
\bjournal{Extremes}
\bvolume{8}
\bpages{237--266}.
\bid{doi={10.1007/s10687-006-0004-0}, issn={1386-1999}, mr={2324891}}
\bptok{imsref}%
\end{barticle}
\endbibitem

\bibitem{surgailis93stable}
\begin{barticle}[mr]
\bauthor{\bsnm{Surgailis},~\bfnm{Donatas}\binits{D.}},
  \bauthor{\bsnm{Rosi{\'n}ski},~\bfnm{Jan}\binits{J.}},
  \bauthor{\bsnm{Mandrekar},~\bfnm{V.}\binits{V.}} \AND
  \bauthor{\bsnm{Cambanis},~\bfnm{Stamatis}\binits{S.}}
(\byear{1993}).
\btitle{Stable mixed moving averages}.
\bjournal{Probab. Theory Related Fields}
\bvolume{97}
\bpages{543--558}.
\bid{doi={10.1007/BF01192963}, issn={0178-8051}, mr={1246979}}
\bptok{imsref}%
\end{barticle}
\endbibitem

\bibitem{takahashi71invariant}
\begin{barticle}[mr]
\bauthor{\bsnm{Takahashi},~\bfnm{Wataru}\binits{W.}}
(\byear{1971}).
\btitle{Invariant functions for amenable semigroups of positive contractions on
  {$L\sp{1}$}}.
\bjournal{K\=odai Math. Sem. Rep.}
\bvolume{23}
\bpages{131--143}.
\bid{issn={0023-2599}, mr={0296256}}
\bptok{imsref}%
\end{barticle}
\endbibitem

\bibitem{tran76problem}
\begin{barticle}[mr]
\bauthor{\bsnm{Tran},~\bfnm{Lanh~Tat}\binits{L.~T.}}
(\byear{1976/77}).
\btitle{On a problem posed by {O}rey and {P}ruitt related to the range of the
  {$N$}-parameter {W}iener process in {$R\sp{d}$}}.
\bjournal{Z. Wahrsch. Verw. Gebiete}
\bvolume{37}
\bpages{27--33}.
\bid{mr={0443063}}
\bptok{imsref}%
\end{barticle}
\endbibitem

\bibitem{wang10association}
\begin{barticle}[mr]
\bauthor{\bsnm{Wang},~\bfnm{Yizao}\binits{Y.}} \AND
  \bauthor{\bsnm{Stoev},~\bfnm{Stilian~A.}\binits{S.~A.}}
(\byear{2010}).
\btitle{On the association of sum- and max-stable processes}.
\bjournal{Statist. Probab. Lett.}
\bvolume{80}
\bpages{480--488}.
\bid{doi={10.1016/j.spl.2009.12.001}, issn={0167-7152}, mr={2593589}}
\bptok{imsref}%
\end{barticle}
\endbibitem

\bibitem{wang10structure}
\begin{barticle}[mr]
\bauthor{\bsnm{Wang},~\bfnm{Yizao}\binits{Y.}} \AND
  \bauthor{\bsnm{Stoev},~\bfnm{Stilian~A.}\binits{S.~A.}}
(\byear{2010}).
\btitle{On the structure and representations of max-stable processes}.
\bjournal{Adv. in Appl. Probab.}
\bvolume{42}
\bpages{855--877}.
\bid{doi={10.1239/aap/1282924066}, issn={0001-8678}, mr={2779562}}
\bptok{imsref}%
\end{barticle}
\endbibitem

\end{thebibliography}
\end{document}